\documentclass[reqno]{amsart}
\usepackage[centertags]{amsmath}
\usepackage{amsfonts}
\usepackage{amssymb}
\usepackage{amsthm}
\usepackage{newlfont}
\usepackage[top=1in, bottom=1in, left=1in, right=1in]{geometry}
\usepackage{mathrsfs}
\usepackage{dsfont}

\theoremstyle{plain}
\newtheorem{ax}{Axiom}
\newtheorem*{ax*}{Axiom}
\newtheorem*{conj}{Conjecture}
\theoremstyle{remark}
\newtheorem*{remark}{Remark}

\newcommand{\Z}{\mathbb{Z}}
\newcommand{\C}{\mathbb{C}}

\DeclareMathOperator{\Fix}{Fix}

\newcommand{\h}[2]{\mathscr{H}_{#1,#2}}
\newcommand{\q}[1]{\mathscr{Q}_{#1}}
\newcommand{\Q}{\mathscr{Q}}
\newcommand{\jac}{\mathcal{J}}

\newcommand{\J}{\mathcal{J}}
\newcommand{\FJRW}[2]{\mathscr{H}_{#1,#2}}

\newcommand{\g}{g}
\newcommand{\lb}{l}

\newcommand{\1}{\mathds{1}}

\newcommand{\set}[1]{\left\{ #1 \right\}}

\newcommand{\br}[1]{\left\langle #1 \right\rangle}

\newcommand{\parlengths}{\setlength{\parindent}{0pt}}
\begin{document}
\date{\today}

\title{FJRW-Rings and Mirror Symmetry}

\author{Marc Krawitz}
\thanks{M.K. is partially Supported by the National Research Foundation of South Africa}
\address{Department of Mathematics, University of Michigan, Ann Arbor, MI , USA}
\email{mkrawitz@umich.edu}

\author{Nathan Priddis}
\address{Department of Mathematics, Brigham Young University, Provo, UT 84602, USA}
\email{priddisn@gmail.com}

\author{Pedro Acosta}
\address{Department of Mathematics, Brigham Young University, Provo, UT 84602, USA}
\email{pedroacosta.lds@gmail.com}

\author{Natalie Bergin}
\address{Department of Mathematics, Brigham Young University, Provo, UT 84602, USA}
\email{natwilde@gmail.com}

\author{Himal Rathnakumara}
\address{Department of Mathematics, Brigham Young University, Provo, UT 84602, USA}
\email{himal46@gmail.com}

\begin{abstract}
The Landau-Ginzburg Mirror Symmetry Conjecture states that for a quasi-homogeneous singularity $W$ and a group $G$ of symmetries of $W$, there is a dual singularity $W^T$ such that the orbifold A-model of $W/G$ is isomorphic to the B-model of $W^T$.  The Landau-Ginzburg A-model is the Frobenius algebra $\h WG$ constructed by Fan, Jarvis, and Ruan, and the B-model is the orbifold Milnor ring of $W^T$.
We verify the Landau-Ginzburg Mirror Symmetry Conjecture for Arnol'd's list of unimodal and bimodal quasi-homogeneous singularities with $G$ the maximal diagonal symmetry group, and include a discussion of eight axioms which facilitate the computation of FJRW-rings.
\end{abstract}

\maketitle\parlengths

\section{Introduction}\label{intro}
In this paper we verify the Landau-Ginzburg Mirror Symmetry
Conjecture for Arnol'd's list of unimodal and bimodal singularities \cite[pg 25]{Arn-1}. Briefly, the conjecture states that for non-degenerate quasi-homogeneous singularities there is a mirror dual singularity such that the ring constructed by Fan-Jarvis-Ruan \cite{FJR} for one is isomorphic to the Landau-Ginzburg B-model (Milnor ring) of the other. The conjecture has already been proven for the simple and parabolic singularities in \cite{FJR}.

The Landau-Ginzburg B-model is an orbifolded Milnor ring.  When the orbifold group is trivial, this is just the classical Milnor ring (local algebra) of the singularity. \cite{Arn-1} 

In \cite{FJR}, Fan, Jarvis, and Ruan construct a cohomological field theory which gives the A-model Frobenius algebra when restricted to genus zero with three marked points. Since the original motivation for the theory was to study generalizations of the Witten equation, we call the A-model Frobenius Algebra the Fan-Jarvis-Ruan-Witten ring, or FJRW-ring, for short.

The singularities for this particular theory are required to be non-degenerate, quasi-homogeneous (i.e. weighted homogeneous) polynomials, having an isolated singularity at the origin. 
Not all of the singularities in the list in \cite{Arn-1} are quasi-homogeneous, so we have only used those which are. Several of the families in \cite{Arn-1} depend on certain parameters, but are only non-degenerate and quasi-homogeneous for particular parameter values.  In these cases, we fix the appropriate parameter values without further comment.

\subsection{Outline of Paper}
\begin{itemize}
\item Section \ref{intro}: Introduction
	\begin{itemize}
	\item \ref{constr} Review of Construction of FJRW-rings
	\item \ref{sum} Additional Notation
	\item \ref{ex} An Example
	\item \ref{form} Format of Results
	\end{itemize}
\item Section \ref{comp}: Computations
	\begin{itemize}
	\item \ref{uni} Exceptional families of unimodal singularities
	\item \ref{bi} Exceptional families of bimodal singularities
	\item \ref{co3} Bi-modal singularities of corank 3
	\item \ref{co2} Bi-modal singularities of corank 2
	\end{itemize}

\end{itemize}

\subsection{Review of Construction} \label{constr}

Let $W$ be a non-degenerate quasi-homogeneous polynomial in the variables $x_1, x_2, \dots, x_N$ with weights $q_1,q_2,\dots,q_N$ respectively. Non-degeneracy requires that these weights are uniquely determined by the condition that each monomial in $W$ has total weight $1$, and that $W$ has an isolated singularity at the origin.

To each quasi-homogeneous polynomial $W$ we can associate a matrix, $B_W$, such that the columns correspond to the variables, and the rows correspond to the terms of the polynomial (see \cite{berglund}). In other words, the entry $(B_W)_{i,j}$ is the power of $x_j$ in the $i$-th monomial of $W$. If the number of monomials coincides with the number of variables, the matrix $B_W$ is square, and the non-degeneracy condition implies that $B_W$ is invertible.  In such cases, we call the polynomial $W$ an \emph{invertible potential}.  Note that the variables of an invertible potential can always be rescaled so the coefficients of the monomials are all equal to one.

To illustrate the correspondence $W\leftrightarrow B_W$, consider the polynomial $W=x^3+xy^4+yz^2$. The corresponding matrix is
\[\left(
\begin{matrix}
3 & 0 & 0\\
1 & 4 & 0\\
0 & 1 & 2
\end{matrix}
\right).
\]

If the matrix $B_W$ is square, its transpose $B_W^T$ will also correspond to a quasi-homogeneous polynomial, which we denote by $W^T$.  In many cases $W^T$ polynomial also has an isolated singularity at the origin, thus satisfying the non-degeneracy condition.

The central charge of $W$ is defined to be
\[
\hat c:=\sum_{j=1}^N(1-2q_j).
\]
The Jacobian ideal $\jac$ is defined by \[\jac=\left(\frac{\partial W}{\partial x_1},\frac{\partial W}{\partial x_2},\dots,\frac{\partial W}{\partial x_N}\right).
\]
The Milnor ring $\q{W}$ is given by
\[
\q{W} :=\C[x_1,x_2,\dots,x_N]/\jac
\]
together with the residue pairing. $\q{W}$ is a finite dimensional vector space over $\C$, with dimension
\[
\mu=\prod_{j=1}^N \left(\frac 1{q_j}-1\right).
\]
It is graded by weighted degree, and the elements of top degree form a one-dimensional subspace generated by $\text{hess}(W)=\det\left(\tfrac{\partial^2 W}{\partial x_i\partial x_j}\right)$.  One can check directly that the top degree is equal to $\hat c$.

For $f,g\in\q{W}$, the residue pairing $\br{f,g}$ may be defined by
\[fg = \frac{\br{f,g}}{\mu}\text{hess}(W) + \text{ lower order terms.}\]

This pairing is non-degenerate, and endows the Milnor ring with the structure of a Frobenius algebra.

To define the FJRW ring, we require in addition to $W$ a choice of a group of diagonal symmetries of $W$. The choice of group heavily affects the resulting structure of the FJRW ring. The maximal group of diagonal symmetries is defined as
\[
G_W=\set{(\alpha_1,\alpha_2,\dots,\alpha_N)\subseteq (\C^*)^N \,|\,W(\alpha_1x_1,\alpha_2x_2,\dots,\alpha_Nx_N)=W(x_1,x_2,\dots,x_N)}
\]
Note that $G_W$ always contains the exponential grading element $J=(e^{2\pi iq_1},e^{2\pi iq_2},\dots,e^{2\pi iq_N})$. In general, the theory requires that the symmetry group be \emph{admissible} (see \cite{FJR} section 2.3). In our computations, we will use either the maximal diagonal symmetry group $G_W$ or the subgroup generated by $J$; both of these groups known to be admissible.

The Landau-Ginzburg Mirror Symmetry Conjecture states the following:
\begin{conj}For a non-degenerate, quasi-homogeneous, invertible singularity $W$ and (maximal) diagonal symmetry group $G$, there is a dual singularity $W^T$ so that the FJRW-ring of $W/G$ is isomorphic to the (unorbifolded) Milnor ring of $W^T$.\end{conj}

\begin{remark}
We use the notation $W^T$ suggestively for the dual singularity, as we identify in this paper a class of examples for which the Berglund-H\"ubsch transposed singularity is the appropriate dual in the context of the Landau-Ginzburg Mirror Symmetry Conjecture.
\end{remark}

We now outline the definition of $\h WG$ as a $\C$-vector space, after which we will define the pairing, grading, and multiplication that make $\h WG$ a Frobenius algebra.

In \cite{FJR}, the state space $\h WG$ is defined in terms of Lefschetz thimbles.  For computational convenience, we give a presentation in terms of Milnor rings, but we should point out that the isomorphism between the two presentations is not canonical.

Let $G$ be an admissible group. For $h\in G$, let $\Fix h\subset\C^N$ be the fixed locus of $h$, and let $N_h$ be its dimension. Define
\[
\mathscr{H}_{h}:=\Omega^{N_h}(\C^{N_h})/\left(dW|_{\Fix h}\wedge \Omega^{N_h-1}\right)\cong\q{W|_{\Fix h}}\cdot\omega
\]
where $\omega=dx_{i_1}\wedge dx_{i_2}\wedge\dots\wedge dx_{i_{N_h}}$ is the natural choice of volume form.

$G$ acts on $\mathscr{H}_{h}$ via its action on the coordinates, and the state space of the FJRW-ring is the vector space of invariants under this action, i.e.
\[
\h WG := \left(\bigoplus_{h\in G}\mathscr{H}_h\right)^G.
\]

$\h WG$ is $\mathbb{Q}$-graded by the so-called $W$-degree, which depends only on the $G$-grading.  To define this grading, first note that each element $h\in G$ can be uniquely expressed as
\[
h=(e^{2\pi i\Theta_1^h},e^{2\pi i\Theta_2^h},\dots,e^{2\pi i\Theta_N^h})
\]
with $0\leq \Theta_i^h < 1$. 

For $\alpha_h \in (\mathscr{H}_{h})^G$, the $W$-degree of $\alpha_h$ is defined by
\begin{equation}
\deg_W(\alpha_h):=N_h+2\sum_{j=1}^N(\Theta_j^h-q_j).\label{degw}
\end{equation}

Since $\Fix h=\Fix h^{-1}$, we have $\mathscr{H}_h \cong \mathscr{H}_{h^{-1}}$, and the pairing on $\q{W|_{\Fix h}}$ induces a pairing
\[(\mathscr{H}_h)^G \otimes (\mathscr{H}_{h^{-1}})^G\to\C.\]
The pairing on $\h WG$ is  the direct sum of these pairings.  Fixing a basis for $\h WG$, we denote the pairing by a matrix $\eta_{\alpha,\beta}=\br{\alpha,\beta}$, with inverse $\eta^{\alpha,\beta}$.

For each pair of non-negative integers $g$ and $n$, with $2g-2+n>0$, the FJRW cohomological field theory produces classes $\Lambda_{g,n}^W(\alpha_1,\alpha_2,\dots,\alpha_n)\in H^*(\overline {\mathscr{M}}_{g,n})$  of complex codimension $D$ for each $n$-tuple $(\alpha_1,\alpha_2,\dots,\alpha_n)\in (\h WG)^n$. The codimension $D$ is given by
\[
D:=\hat c_W(g-1)+\frac 12\sum_{i=1}^n\deg_W(\alpha_i),
\]
and the $n$-point correlators are defined to be
\[
\br{\alpha_1,\dotsc,\alpha_n}_{g,n}:=\int_{\overline {\mathscr{M}}_{g,n}}\Lambda_{g,n}^W(\alpha_1,\dotsc,\alpha_n).
\]
The correlator $\br{\alpha_1,\dotsc,\alpha_n}_{g,n}$ vanishes unless the codimension of $\Lambda_{g,n}^W(\alpha_1,\dotsc,\alpha_n)$ is zero. The ring structure on $\mathscr{H}_{W,G}$ is determined by the genus-zero three-point correlators. In other words, if $r,s\in \h{W}{G}$, then
\begin{equation}
r*s:=\sum_{\alpha,\beta}\br{r,s,\alpha}_{0,3}\eta^{\alpha,\beta}\beta \label{mult}
\end{equation}
where the sum is taken over all choices of $\alpha$ and $\beta$ in a fixed basis of $\h WG$.

The classes $\Lambda_{g,n}^W(\alpha_1,\dotsc,\alpha_n)$ satisfy the following axioms that allow us to compute most of the three-point correlators $\br{\alpha_1,\alpha_2,\alpha_3}$ explicitly.

\begin{ax} Dimension: If $D\notin \frac 12 \Z$, then $\Lambda_{g,n}^W(\alpha_1,\alpha_2,\dots,\alpha_n)=0$. Otherwise, $D$ is the complex codimension of the class $\Lambda_{g,n}^W(\alpha_1,\alpha_2,\dots,\alpha_n)$. In particular, if $g=0$ and $n=3$, then $\br{\alpha_1,\alpha_2,\alpha_3}=0$ unless $D=0$.
\end{ax}

Notice that in the case where $g=0$ and $n=3$, $D=0$ if and only if $\sum_{i=1}^3\deg_W\alpha_i=2\hat c$.

\begin{ax} Symmetry:
Let $\sigma\in S_n$. Then
\[
\br{\alpha_1,\dotsc,\alpha_n}_{g,n}=\br{\alpha_{\sigma(1)},\dotsc,\alpha_{\sigma(n)}}_{g,n}.
\]
\end{ax}

The next few axioms rely on the degrees of line bundles $\mathscr{L}_1,\dotsc,\mathscr{L}_N$ endowing an orbicurve with a so-called \emph{$W$-structure}; however, this can be reduced to a simple numerical criterion. Consider the class $\Lambda_{g,k}^W(\alpha_1,\alpha_2,\dots,\alpha_k)$, with $\alpha_j\in(\mathscr{H}_{h_j})^G$ for each $j\in\{1,\dots,N\}$. For each variable $x_j$, define $\lb_j$ by
\[
\lb_j=q_j(2g-2+k)-\sum_{i=1}^k\Theta_j^{h_i}
\]

\begin{ax} Integer degrees: If $\lb_j\notin \Z$ for some $j\in\set{1,\dots,N}$, then $\Lambda_{g,k}^W(\alpha_1,\alpha_2,\dots,\alpha_k)=~0$.
\end{ax}

\begin{ax}
Concavity: If $\lb_{j}<0$ for all $j\in\set{1,2,3}$, then $\br{\alpha_1,\alpha_2,\alpha_3}=1$.
\end{ax}

The next axiom is related to the Witten map:
\begin{eqnarray*}
\mathcal W:\bigoplus_{\substack{j=1\\}}^N \C^{h^0_j}\rightarrow \bigoplus_{j=1}^N \C^{h^1_j}\\
\mathcal W=\left(\overline{\frac{{\partial W}}{{\partial x_1}}}, \overline{\frac{{\partial W}}{{\partial x_2}}}, \dots, \overline{\frac{{\partial W}}{{\partial x_N}}}\right)
\end{eqnarray*}

where $h^0_j$ and $h^1_j$ are defined by
\[
h^0_j:=\begin{cases}
0       & \text{if } \lb_j< 0\\
\lb_j+1 & \text{if }\lb_j \geq 0
\end{cases}
\]
\[
h^1_j:=\begin{cases}
-\lb_j-1 & \text{if } \lb_j<0\\
0        & \text{if } \lb_j\geq 0
\end{cases}
\]
so that both are non-negative integers satisfying $h^0_j - h^1_j = l_j + 1$.\footnote{The reader may note that $h_j^i$ is just the dimension of the $i$-th cohomology of the $j$-th line bundle in the $W$-structure, and this relation is just Riemann-Roch for a line-bundle $L_j$ of  degree $l_j$ on $\C P^1$.  The preceding definitions simply serve to axiomatize the entire construction.}

The fact that the Witten map is well-defined is a consequence of the geometric conditions on the $\mathscr{L}_j$ considered in \cite{FJR}.  For further details, we refer readers to the original paper.

If $\Lambda_{g,n}^W(\alpha_1,\dots,\alpha_n)$ is a class of codimension zero, we obtain a complex number by integrating over $\overline{\mathscr{M}}_{g,n}$.  Abusing notation, we will refer to the class $\Lambda_{g,n}^W(\alpha_1,\dots,\alpha_n)$ and its integral over $\overline{\mathscr{M}}_{g,n}$ interchangeably.

\begin{ax}
Index Zero: Consider the class $\Lambda_{g,n}^W(\alpha_1,\alpha_2,\dots,\alpha_n)$, with $\alpha_i\in \h{\gamma_i}{G}$. If $\Fix \gamma_i=\set 0$ for each $i\in\set{1,2,\dots,n}$ and 
\[
\sum_{j=1}^N (h^0_j - h^1_j) = 0
\]
then $\Lambda_{g,n}(\alpha_1,\alpha_2,\dots,\alpha_n)$ is of codimension zero and $\Lambda_{g,n}^W(\alpha_1,\alpha_2,\dots,\alpha_n)$ is equal to the degree of the Witten map. 
\end{ax}

\begin{ax}
Composition: If the four-point class, $\Lambda_{g,n}^W(\alpha_1,\alpha_2,\alpha_3,\alpha_4)$ is of codimension zero, then it decomposes in terms of three-point correlators in the following way:
\[
\Lambda_{0,4}^W(\alpha_1,\alpha_2,\alpha_3,\alpha_4)=\sum_{\beta,\delta}\br{\alpha_1, \alpha_2,\beta}\eta^{\beta,\delta}\br{\delta,\alpha_3,\alpha_4}.\]
\end{ax}

Note that $\Fix J=\set 0$ so $\mathscr{H}_{J}\cong\C$. The identity element in the FJRW-ring is an element of $\mathscr{H}_{J}$, and we denote this element by $\1$.

\begin{ax}
Pairing: For $\alpha_1,\alpha_2\in \h WG$,  $\br{\alpha_1,\alpha_2,\mathds{1}}=\eta(\alpha_1,\alpha_2)$, where $\eta$ is the pairing in $\mathscr{H}_{W,G}$.
\end{ax}

\begin{ax} \label{ax:sums} Sums of singularities: If $W_1\in\C[x_1,\dots,x_r]$ and $W_2\in\C[y_1,\dots,y_s]$ are two non-degenerate, quasi-homogeneous polynomials with maximal symmetry groups $G_1$ and $G_2$, then the maximal symmetry group of $W=W_1+W_2$ is $G=G_1\times G_2$, and there is an isomorphism of Frobenius algebras
\[
\h{W}{G}\cong \h{W_1}{G_{W_1}}\otimes \h{W_2}{G_{W_2}}
\]
\end{ax}

\begin{remark}\label{rem:sums}
We note an important consequence of Axiom \ref{ax:sums}.  Under the same hypotheses as in the statement of the axiom, we have a Frobenius Algebra isomorphism
\[\Q_{W}\cong \Q_{W_1}\otimes \Q_{W_2},\]
and similarly
\[\Q_{W^T}\cong \Q_{W_1^T}\otimes \Q_{W_2^T}.\]
Consequently, in order to prove the Mirror Symmetry Conjecture for $W=W_1+W_2$, a sum of decoupled polynomials (with maximal $A$-model orbifold group), it suffices to prove it for $W_1$ and $W_2$ individually.
\end{remark}

These axioms allow us to compute most of the three-point correlators of the FJRW-rings. In some cases, the axioms are not enough to compute all of the correlators; however, in most of these cases, one can still verify the mirror symmetry conjecture.

\subsection{Additional Notation}\label{sum}
A singularity is said to be \emph{invertible} if the number of monomials equals the number of variables.
We have found that the Landau-Ginzburg mirror symmetry conjecture holds for the invertible unimodal and bimodal singularities orbifolded by the maximal group of diagonal symmetries, $G_W$. The conjecture was verified in \cite{FJR} for the simple singularities.  

In most cases considered, the maximal symmetry group, $G_W$, is cyclic. In these cases we have adopted the following notation. Let $\g$ be a generator for $G_W$. If $J$ generates $G_W$, take $g=J$.
If $\Fix{\g^k}=\set 0$, define
\[
e_k=1\in \mathscr{H}_{\g^k} \cong \C,
\]
otherwise if $\Fix{\g^k}= \C x_{i_1}\oplus\dotsb\oplus\C x_{i_{N_g}}$ define
\[
e_k=dx_{i_1}\wedge dx_{i_2}\wedge\dots\wedge dx_{i_{N_g}}\in \mathscr{H}_{\g^k}.
\]
We denote coordinate subspaces of $\C^N$ with subscripts indicating the non-zero variables in these subspaces.  So, for example, the $xy$-plane in $\C^3$ will be denoted by $\C^2_{xy}$.

Our computations are often made easier by judicious use of associativity. This will be reflected typographically in the grouping of terms. So, for example, 
\[
\alpha*\alpha\beta
\]
indicates that $\alpha\beta$ should be computed first, and then multiplied by $\alpha$. 

\subsection{An example}\label{ex}
We will now give an example demonstrating the construction and our methods of computation more fully.

\subsubsection{$E_{19}$ with maximal symmetry group}

The polynomial for $E_{19}$ is $x^3+xy^7$.
The corresponding weights for each variable are $q_x=\frac 13, q_y=\frac 2{21}$ and the central charge is $\hat c=\frac {24}{21}$.
The Jacobian ideal is $\jac=(3x^2 + y^7,7xy^6)$.

The maximal group of diagonal symmetries is given by
\[
G=\br{(\alpha,\beta)\mid \alpha^3, \alpha\beta^7}
\]
From the relations, we can see that $\alpha=\beta^{-7}$, so that $|G|=21$ and $G$ is cyclic. The exponential grading element is $J=(e^{2\pi i\tfrac{1}{3}}, e^{2\pi i\tfrac{2}{21}})$, which has order 21, so the maximal group of diagonal symmetries is generated by $J$.

The fixed point locus of $J^k$ is given by
\[\Fix{J^k}=\begin{cases}
\C^2 & \text{if}\quad k=0\\
\C_x & \text{if}\quad 3|k, k\neq 0\\
0    & \text{otherwise}
\end{cases}
\]

The Milnor ring is $\C[x,y]/\jac\cong \br{1,x,x^2,y,\dots,y^{6},xy,xy^2,\dots,xy^5,x^2y,x^2y^2,\dots,x^2y^5}$ for $J^0=I$. If $3|k$, then $W\mid_{\Fix J^k}=x^3$. And in the case that $3\nmid k$, the Milnor ring is trivial. So the Milnor rings are given by
\[
\q{}|_{\Fix J^k}=
\begin{cases}
\br{1,x,x^2,y,\dots,y^{6},xy,xy^2,\dots,xy^5,x^2y,x^2y^2,\dots,x^2y^5}, \mu=19 & \text{if}\quad k=0\\
\br{1,x}, \mu=2 & \text{if}\quad 3|k, k\neq 0\\
\br 1  & \text{otherwise}
\end{cases}
\]

Now for each choice of $k$ with $3\nmid k$, $\Fix J^k = \{0\}$ so $\h{J^k}G\cong\C$ with trivial $G$-action. So $e_k$ is a generator for the state space. To take $G$-invariants, we need only consider the action of $J$, since $G=\langle J \rangle$.

If $k\neq 0$ and $3\mid k$, then we must consider the action of $J$ on $x^j\,dx$ for $j\in \{0,1\}$.
\[
J\cdot x^j\,dx=\exp\left(2\pi i\tfrac{(j+1)}{3}\right)x^j\,dx
\]
This action clearly has no invariants.

For $k=0$, we must consider the action of $J$ on the Milnor ring associated to $J^0$.
\[
J\cdot x^ly^m\,dx\wedge dy=\exp\left(2\pi i\tfrac{(7l+7+2m+2)}{21}\right) x^ly^m\,dx\wedge dy.
\]
This action has a single invariant element, namely $y^6\,dx\wedge dy=y^6e_0$. So a basis for the state space is
\[
\h{E_{19}}{\br
J}=\br{y^6e_{0},e_1,e_2,e_4,e_5,e_7,e_8,e_{10},e_{11},e_{13},e_{14},e_{16},e_{17},e_{19},e_{20}}.
\]
\vspace{\baselineskip}

The $W$-degrees are straightforward to compute from formula (\ref{degw}) on pg 3. Again, notice that the $W$-degree only depends on the power of the generator $J$. We give the invariants of each sector in the following table as well.

\begin{center}
\begin{tabular}{|c||c|c|c|c|c|c|c|c|c|c|c|c|c|c|c|c|c|c|c|c|c|c|c|c|}
\hline
$k$        & 0 & 1 & 2 & 4 & 5 & 7 & 8 & 10& 11 & 13 & 14 & 16 & 17 & 19 & 20 \\
\hline
$|G|\cdot\deg_W$   & 24 & 0 & 18 & 12 & 30 & 24 & 42 & 36& 12 & 6 & 24 & 18 & 36 & 30 & 48\\
\hline
invariants &$y^6e_{0}$& $\1$& $e_2$& $e_4$& $e_5$& $e_7$& $e_8$& $e_{10}$& $e_{11}$& $e_{13}$& $e_{14}$& $e_{16}$& $e_{17}$& $e_{19}$& $e_{20}$\\
\hline
\end{tabular}

\end{center}

\vspace*{\baselineskip}


Several three-point correlators can be computed using the pairing axiom. Note that all twisted sectors correspond to trivial fixed loci, so the residue pairing between twisted sectors is given by
\begin{gather*}
\mathscr{H}_{J^k}\otimes\mathscr{H}_{J^{21-k}} \stackrel{\br{\,,\,}}{\longrightarrow} \C\\
\br{e_k,\, e_{21-k} } = 1
\end{gather*}

So, by the pairing axiom, $\br{\1,\1,e_{20}}$, $\br{\1,e_2,e_{19}}$,
$\br{\1,e_4,e_{17}}$, $\br{\1,e_5,e_{16}}$, $\br{\1,e_7,e_{14}}$,
$\br{\1,e_8,e_{13}}$, and $\br{\1,e_{10},e_{11}}$ are all equal to 1.

The pairing in the untwisted sector gives us
\[
\br{y^6e_0,y^6e_0,\1}=\eta_{y^6e_0, y^6e_0}=\frac{19y^{12}}{\mbox{Hess}}=\frac{19y^{12}}{-133y^{12}}=-1/7.
\]

The following three-point correlators have $\lb_x$ and $\lb_y$ both less than zero: $\br{e_2,e_4,e_{16}}$, $\br{e_2,e_7,e_{13}}$, $\br{e_4,e_4,e_{14}}$, $\br{e_4,e_5,e_{13}}$, $\br{e_4,e_7,e_{11}}$, $\br{e_{11},e_{13},e_{19}}$, $\br{e_{11},e_{16},e_{16}}$, $\br{e_{13},e_{13},e_{17}}$, $\br{e_{13},e_{14},e_{16}}$. Therefore, by the concavity axiom each is equal to 1.

For the three-point correlator $\br{y^6e_0,e_{11},e_{11}}$, we apply the composition axiom to the class $\Lambda^W_{0,4}(e_{11},e_{11},e_{11},e_{11})$. The values for the line bundle degrees are $\lb_x=-2$ and $\lb_y=0$. So we have $H^0=0\oplus \C_y$ and $H^1=\C_x\oplus 0$. The Witten map $H^0\to H^1$ is given by the complex conjugate of the gradient of $W$. In other words it maps
\[(0,y)\mapsto (\overline y^7, 0).\]
The degree of this map is $-7$, so the composition axiom tells us that
\[
-7=\Lambda_{0,4}^W(e_{11},e_{11},e_{11},e_{11})=\sum_{\alpha,\beta}\br{e_{11}, e_{11},\alpha}\eta^{\alpha,\beta}\br{\beta,e_{11},e_{11}},
\]
where $\alpha$ and $\beta$ range over a fixed basis for the FJRW-ring.

By degree considerations, the only non-zero contribution to this sum  occurs when $\alpha = y^6e_0 = \beta$.  Setting $a:=\br{y^6e_0,e_{11},e_{11}}$, we get the following equation:
\begin{align}\label{e11}
{-7} &= \br{e_{11}, e_{11},y^6e_0}\eta^{y^6e_0,y^6e_0}\br{y^6e_0,e_{11},e_{11}}\\
   &= {-7} a^2\notag,
\end{align}
where the ${-7}$ on the right is the contribution to the inverse of the pairing corresponding to $\br{y^6e_0,\,y^6e_0,\1}={-1/7}$ computed above.

From this we can see $\br{y^6e_0,e_{11},e_{11}}=\pm 1$. We will see below that either choice gives a Frobenius algebra that is isomorphic to $\q{E_{19}^T}$, where $E_{19}^T = x^3y+y^7$.

This completes the computation of the three-point correlators. All others are required to vanish by Axioms 1 and 2. 

\vspace*{\baselineskip}
Recall that mulplication in $\h{E_{19}}{G}$ is given by equation (\ref{mult}).
Examining degrees we see that $e_{13}$ is a generator for $\h{E_{19}}{G}$. We compute
\begin{align*}
e_{13}^2 &=\sum{\alpha,\beta}\br{e_{13},e_{13},\alpha}\eta^{\alpha,\beta}\beta\\
         &=\br{e_{13},e_{13},e_{17}}\eta^{e_{17},e_4}e_4\\
         &= e_4
\end{align*}
where the first summation is taken over all $\alpha$ and $\beta$ among the generators we have listed for the state space. The second equality follows because the only non-zero three-point correlator with $e_{13}$ occurring twice is $\br{e_{13},e_{13},e_{17}}=1$. Again examining degrees, we see that $e_{11}$ is also a generator of $\h{E_{19}}{G}$.

One can check directly that $e_{13}$ and $e_{11}$ generate $\h{E_{19}}{G}$ as a Frobenius algebra, so we may define a surjective map $\varphi:\C[X,Y]\rightarrow \h{E_{19}}{G}$ by $X\mapsto e_{11}$ and $Y\mapsto e_{13}$, and extend to a $\C$-algebra homomorphism.

We see that $e_{11}^2*e_{13}= 0$, $e_{11}^3= -7e_{10}$, and $e_{13}^6=e_{10}$. Hence $(3X^2Y,X^3+7Y^6)\subset \ker\varphi$.

 Therefore since $\C[X,Y]/(3X^2Y,X^3+7Y^6)=\q{E_{19}^T}$ has dimension 15, the same as $\h{E_{19}}{G}$, we deduce that the inclusion $(3X^2Y,X^3+7Y^6)\subset \ker\varphi$ is an equality, therefore the map induces a degree-preserving isomorphism $\q{E_{19}^T}\cong \h{E_{19}}{G}$.

\subsection{Format of results}\label{form}

For each singularity, we will display the information in the following pattern:
\begin{itemize}

\item Name of singularity, the defining polynomial, the Jacobian ideal, the weights associated to each variable, and the central charge. We will also give the symmetry group used in the construction, typically $G_W$, which we will henceforth denote by $G$.

\item Fixed locus for each group element.

\item Basis for the Milnor ring of $W$ restricted to each fixed locus. 

\item Table of sectors with non-trivial $G$-invariants, including the invariant elements and their $W$-degrees.  For clarity of exposition, we will multiply $W$-degrees by a factor of $|G|$.

\item We will give the values of the three-point correlators that are not required to vanish by Axioms 1 and 2. These will be grouped in the following order: those computed by the Pairing axiom, those by the Concavity axiom, those by the Index Zero axiom, and those by the Composition axiom.  Any correlators for which the axioms do not suffice will be listed last, including any relations among them.

\item Finally, we will describe, where possible, an isomorphism between the FJRW-ring of $W$ and the Milnor ring of $W^T$.
\end{itemize}

\section{Computations}\label{comp}

We take our examples from the unimodal and bimodal singularities listed by Arnol'd \cite{Arn-1}. Many of these singularities are quasi-homogeneous only after fixing specific parameter values, which we do without further comment.

\subsection{Unimodal Singularities}\label{uni}

\subsection*{$\mathbf{Q_{10} = x^2z + y^3 + z^4}$.}
Axiom \ref{ax:sums} applies here.  By the subsequent comment, it suffices to prove the Mirror Symmetry Conjecture for $D_5 = x^2z+z^4$ and $A_2 = y^3$.  The conjecture was proved for the simple (ADE) singularities in \cite{FJR}, so it holds in this case also.
\newline\hrule

\subsection*{$\mathbf{Q_{11} = x^2z + y^3 + yz^3}$.}
\[\J = (2xz,\,3y^2+z^3,\,x^2+3yz^2)
\phantom{XXX}q_x =\tfrac{7}{18},\,q_y=\tfrac{6}{18},\, q_z=\tfrac{4}{18},\,\hat c=\tfrac{20}{18}
\phantom{XXX}G = \br{J} = \Z/18\Z\]
\[\Fix J^k = \begin{cases}
  \C^3 & \text{if } k=0\\
  \C_{y} & \text{if } k=3,\,6,\,12,\,15\\
  \C^2_{yz} & \text{if } k=9\\
  0 & \text{otherwise}
\end{cases} 
\phantom{XXX}
\Q|_{\Fix J^k} = \begin{cases}
  \br{1,\,z,\,y,\,x,\,z^2,\,yz,\,z^3,\,xy,\,yz^2,\,z^4,\,z^5},\,\mu=11\\
  \br{1,\,y,\,y^2},\,\mu=3\\
  \br{1,\,z,\,y,\,z^2,\,yz,\,z^3,\,z^4},\,\mu=7\\
  \br{1}\\
\end{cases} \]

\begin{center}
\begin{tabular}{|c||c|c|c|c|c|c|c|c|c|c|c|c|c|}\hline
$k$	&1	&2	&4	&5	&7	&8	&9	&10	&11	&13	&14	&16	&17\\\hline
$|G|\cdot\deg_W$	&0	&34	&30	&28	&24	&22	&20	&18	&16	&12	&10	&6	&40\\\hline
invariants	& $\1$	& $e_{2}$ & $e_{4}$	& $e_{5}$	& $e_{7}$	& $e_{8}$	& $z^2e_{9}$& $e_{10}$	& $e_{11}$	& $e_{13}$	& $e_{14}$	& $e_{16}$	& $e_{17}$\\\hline
\end{tabular}
\end{center}

Potential non-zero correlators:

By the pairing axiom,
$\br{\1, \1, e_{17}}$,
$\br{\1,	e_2, e_{16}}$,
$\br{\1,	e_4, e_{14}}$,
$\br{\1,	e_5, e_{13}}$,
$\br{\1,	e_{7},	e_{11}}$, and
$\br{\1,	e_{8},	e_{10}}$ are all equal to $1$, and
$\br{\1,	z^2e_{9},	z^2e_{9}}=-\tfrac{1}{3}$.

By the concavity axiom,
$\br{e_{5},	e_{16},	e_{16}}$,
$\br{e_{7},	e_{14},	e_{16}}$,
$\br{e_{10},	e_{11},	e_{16}}$,
$\br{e_{10},	e_{13},	e_{14}}$,
are all equal to $1$.

By the index-zero axiom,
$\br{e_{11}, e_{13}, e_{13}}$ and $\br{e_{8}, e_{13}, e_{16}}$ are $-2$.

By the composition axiom,
\begin{equation}
-3 = \br{z^2e_{9},	e_{14},	e_{14}}\eta^{z^2e_9,z^2e_9}\br{z^2e_{9},	e_{14},	e_{14}},\tag{\dag}\label{eq:composition relation}
\end{equation}
where $\eta^{z^2e_9,z^2e_9}=-3$, so $\br{z^2e_{9},	e_{14},	e_{14}}=\pm 1$.

Consider the map
\[\varphi:\C[X,Y,Z]\to \FJRW{Q_{11}}{G}\]
defined by $X\mapsto e_{10}$, $Y\mapsto e_{14}$, and $Z\mapsto e_{16}$ extending by $\C$-linearity and multiplicativity.  One can check directly that this map is surjective.

Our correlators tell us that $e_{16}^3 = -2e_{10}$, $e_{14}^3 = -3 e_{4} = -3 e_{10}*e_{16}^2$, and $e_{14}^2*e_{16}=0$.

Hence $(2X+Z^3,\,3Y^2Z,\,Y^3+3XZ^2)\subseteq \ker \varphi$.  Since
$\C[X,Y,Z]/(2X+Z^3,\,3Y^2Z,\,Y^3+3XZ^2)=\q{Q_{11}^T}$ has dimension 13, we deduce the inclusion is in fact equality, and $\FJRW{Q_{11}}{G}\cong\q{Q_{11}^T}$.
\newline\hrule

\subsection*{$\mathbf{Q_{12} = x^2z + y^3 + z^5}$.}
Axiom \ref{ax:sums} is applicable, and it suffices to prove the Mirror Symmetry Conjecture for $D_6 = x^2z+z^5$ and $A_2 = y^3$.  This follows from the study of simple singularities in \cite{FJR}.
\newline\hrule

\subsection*{$\mathbf{S_{11} = x^2z + yz^2 + y^4}$.}
\[\J = (2xz,\,z^2+4y^3,\,x^2+2yz)
\phantom{XXX} q_x =\tfrac{5}{16},\,q_y=\tfrac{4}{16},\,
q_z=\tfrac{6}{16},\,\hat c=\tfrac{18}{16} \phantom{XXX} G = \br{J} = \Z/16\Z\]
\[\Fix J^k = \begin{cases}
  \C^3_{xyz} & \text{if } k=0\\
  \C_{y} & \text{if } k=4,\,12\\
  \C^2_{yz} & \text{if } k=8\\
  0 & \text{otherwise}
\end{cases}
\phantom{XXX}
\Q|_{\Fix J^k} = \begin{cases}
  \br{1,\,y,\,x,\,z,\,y^2,\,xy,\,x^2,\,z^2,\,xy^2,\,y^2z,\,z^3},\,\mu=11\\
  \br{1,\,z,\,z^2},\,\mu=3\\
  \br{1,\,z,\,y,\,y^2,\,y^3},\,\mu=5\\
  \br{1}\\
\end{cases} \]

\begin{center}
\begin{tabular}{|c||c|c|c|c|c|c|c|c|c|c|c|c|c|}\hline
$k$	&1	&2	&3	&5	&6	&7	&8	&9	&10	&11	&13	&14	&15\\\hline
$|G|\cdot\deg_W$	&0	&30	&28	&24	&22	&20	&18	&16	&14	&12	&8	&6	&36\\\hline
invariants	& $\1$ & $e_{2}	$& $e_{3}$ & $e_{5}$ & $e_{6}$ & $e_{7}$ & $ze_{8}$ & $e_{9}$ & $e_{10}$ & $e_{11}$ & $e_{13}$ & $e_{14}$ & $e_{15}$\\\hline
\end{tabular}
\end{center}

Potential non-zero correlators:

By the pairing axiom,
$\br{\1, \1, e_{15}}$,
$\br{\1,    e_2, e_{14}}$,
$\br{\1,    e_3, e_{13}}$,
$\br{\1,    e_5, e_{11}}$,
$\br{\1,    e_{6},  e_{10}}$, and
$\br{\1,    e_{7},  e_{9}}$ are all equal to $1$, and
$\br{\1,    ze_{8}, ze_{8}}=-\tfrac{1}{2}$.

By the concavity axiom,
$\br{e_{5}, e_{14}, e_{14}}$,
$\br{e_{6}, e_{13}, e_{14}}$,
$\br{e_{9}, e_{10}, e_{14}}$, and
$\br{e_{9}, e_{11}, e_{13}}$
are all equal to $1$.

By the index-zero axiom,
$\br{e_{7}, e_{13}, e_{13}}$,
$\br{e_{10}, e_{10}, e_{13}}$, and
$\br{e_{11}, e_{11}, e_{11}}$
are all equal to ${-2}$.

The correlator $\br{e_{8},  e_{11}, xe_{14}}$ may be non-zero, and the composition axiom can be used to compute it.  However, using associativity of the product we can compute the ring structure on $\FJRW{S_{11}}{G}$ without it.

Consider the map
\[\varphi:\C[X,Y,Z]\to \FJRW{S_{11}}{G}\]
defined by $X\mapsto e_{9}$, $Y\mapsto e_{13}$, and $Z\to e_{14}$ extending by $\C$-linearity and multiplicativity.  One can check directly that this map is surjective.

Our correlators tell us that $e_{13}^2 = -2e_{9}$, $e_{14}^4 = -2e_{13}e_{9}$, and $e_{13}e_{14}*e_{14}^2=0$.

Hence $(2X+Y^2,\, 2XY+Z^4,\, YZ^3)\subseteq \ker \varphi$. Since
$\C[X,Y,Z]/(2X+Y^2,\, 2XY+Z^4,\, YZ^3)=\q{S_{11}^T}$ has dimension 13, we deduce the inclusion is in fact equality, and $\FJRW{S_{11}}{G}\cong\q{S_{11}^T}$.
\newline\hrule

\subsection*{$\mathbf{S_{12} = x^2z + yz^2 + xy^3}$.}
\[\J = (2xz,\,z^2+3xy^2,\,x^2+2yz)
\phantom{XXX} q_x =\tfrac{4}{13},\,q_y=\tfrac{3}{13},\, q_z=\tfrac{5}{13},\,\hat c=\tfrac{15}{13}
\phantom{XXX} G = \br{J} = \Z/13\Z\]
\[\Fix J^k = \begin{cases}
  \C^3_{xyz} & \text{if } k=0\\
  0 & \text{otherwise}
\end{cases}
\phantom{XXX}
\Q|_{\Fix J^k} = \begin{cases}
  \br{1,\,y,\,x,\,z,\,y^2,\,xy,\,yz,\,xz,\,xy^2,\,y^2z,\,xyz,\,xy^2z},\,\mu=12\\
  \br{1} 
\end{cases} \]

\begin{center}
\begin{tabular}{|c||c|c|c|c|c|c|c|c|c|c|c|c|}\hline
$k$	&1	&2	&3	&4	&5	&6	&7	&8	&9	&10	&11	&12\\\hline
$|G|\cdot\deg_W$	&0	&24	&22	&20	&18	&16	&14	&12	&10	&8	&6	&30\\\hline
invariants	& $\1$	& $e_{2}$	& $e_{3}$	& $e_{4}$	& $e_{5}$	& $e_{6}$	& $e_{7}$	& $e_{8}$	& $e_{9}$	& $e_{10}$	& $e_{11}$	& $e_{12}$\\\hline
\end{tabular}
\end{center}

Potential non-zero correlators:

By the pairing axiom,
$\br{\1, \1, e_{12}}$,
$\br{\1,	e_2, e_{11}}$,
$\br{\1,	e_3, e_{10}}$,
$\br{\1,	e_4, e_{9}}$,
$\br{\1,	e_{5},	e_{8}}$, and
$\br{\1,	e_{6},	e_{7}}$ are all equal to $1$.

By the concavity axiom,
$\br{e_{5},	e_{11},	e_{11}}$,
$\br{e_{6},	e_{10},	e_{11}}$,
$\br{e_{7},	e_{9},	e_{11}}$, and
$\br{e_{8},	e_{9},	e_{10}}$
are all equal to $1$.

By the index-zero axiom,
$\br{e_{7}, e_{10}, e_{10}}={-2}$,
$\br{e_{8}, e_{11}, e_{11}}={-2}$, and
$\br{e_{9}, e_{9}, e_{9}}=-{3}$.

Consider the map
\[\varphi:\C[X,Y,Z]\to \FJRW{S_{12}}{G}\]
defined by $X\mapsto e_{9}$, $Y\mapsto e_{10}$, and $Z\to e_{11}$ extending by $\C$-linearity and multiplicativity.  One can check directly that this map is surjective.

Our correlators tell us that $e_{9}^2 = -3e_{4}= -3e_{10}e_{11}^2$, $e_{10}^2 =-2e_6 = -2e_9e_{11}$, and $e_{11}^3 = -2e_{5} = -2e_9e_{10}$.

Hence $(Y^2+2XZ,\,Z^3+2XY,\,X^2+3YZ^2)\subseteq \ker \varphi$. Since
$\C[X,Y,Z]/(Y^2+2XZ,\,Z^3+2XY,\,X^2+3YZ)=\q{S_{12}^T}$ has dimension 13, we deduce the inclusion is in fact equality, and $\FJRW{S_{12}}{G}\cong\q{S_{12}^T}$.
\newline\hrule

\subsection*{$\mathbf{U_{12} = x^3 + y^3 + z^4}$.}
By Axiom 8, the Mirror Symmetry Conjecture holds here because it holds for the simple singularities $A_2 = x^3$ and $A_3 = z^4$.\newline\hrule

\subsection*{$\mathbf{Z_{11} = x^3y + y^5}$.}
\[\J = (3x^2y,\,x^3+5y^4)
\phantom{XXX} q_x =\tfrac{4}{15},\,q_y=\tfrac{3}{15},\,\hat c=\tfrac{16}{15}
\phantom{XXX} G = \br{J} = \Z/15\Z\]
\[\Fix J^k = \begin{cases}
  \C^2_{xy} & \text{if } k=0\\
  \C_{y} & \text{if } k=5,\,10\\
  0 & \text{otherwise}
\end{cases}
\phantom{XXX}
\Q|_{\Fix J^k} = \begin{cases}
  \br{1,\,y,\,x,\,y^2,\,xy,\,x^2,\,y^3,\,xy^2,\,x^3,\,xy^3,\,x^4},\,\mu=11\\
  \br{1,\,y,\,y^2,\,y^3},\,\mu=4\\
  \br{1} 
\end{cases} \]

\begin{center}
\begin{tabular}{|c||c|c|c|c|c|c|c|c|c|c|c|c|c|}\hline
$k$	&0	&1	&2	&3	&4	&6	&7	&8	&9	&11	&12	&13	&14\\\hline
$|G|\cdot\deg_W$	&16	&0	&14	&28	&12	&10	&24	&8	&22	&20	&4	&18	&32\\\hline
invariants	& $x^2e_{0}$ & $\1$	& $e_{2}$	& $e_{3}$	& $e_{4}$	& $e_{6}$	& $e_{7}$	& $e_{8}$	& $e_{9}$	&  $e_{11}$	& $e_{12}$	& $e_{13}$ & $e_{14}$ \\\hline
\end{tabular}
\end{center}

Potential non-zero correlators:

By the pairing axiom,
$\br{\1,	\1, e_{14}}$,
$\br{\1,	e_2, e_{13}}$,
$\br{\1,	e_3, e_{12}}$,
$\br{\1,	e_{4},	e_{11}}$,
$\br{\1,	e_{6},	e_{9}}$, and
$\br{\1,	e_{7},	e_{8}}$ are all equal to $1$, and
$\br{\1, x^2e_{0}, x^2e_{0}}={-\tfrac{1}{3}}$.

By the concavity axiom,
$\br{e_{2},	e_{2},	e_{12}}$,
$\br{e_{2},	e_{6},	e_{8}}$,
$\br{e_{4},	e_{6},	e_{6}}$,
$\br{e_{6},	e_{12},	e_{13}}$,
$\br{e_{7},	e_{12},	e_{12}}$, and
$\br{e_{8},	e_{11},	e_{12}}$
are all equal to $1$.

By the index-zero axiom,
$\br{e_{4}, e_{4}, e_{8}}={-3}$.

The composition axiom can be used to compute
$\br{x^2e_{0},	e_{4},	e_{12}}$, and
$\br{x^2e_{0},	e_{8},	e_{8}}$, but we do not need these to establish the desired isomorphism.

Consider the map
\[\varphi:\C[X,Y]\to \FJRW{Z_{11}}{G}\]
defined by $X\mapsto e_{6}$, and $Y\mapsto e_{12}$ extending by $\C$-linearity and multiplicativity.  One can check directly that this map is surjective.

Our correlators tell us that $e_{12}^3*e_{12}^2 = -3e_{11} = -3e_{6}^2$ and $e_{6}e_{12}^4 = 0$.

Hence $(3X^2+Y^5,\,5XY^4)\subseteq \ker \varphi$.  Since
$\C[X,Y]/(3X^2+Y^5,\,5XY^4)=\q{Z_{11}^T}$ has dimension 13, we deduce the inclusion is in fact equality, and $\FJRW{Z_{11}}{G}\cong\q{Z_{11}^T}$.
\newline\hrule

\subsection*{$\mathbf{Z_{12} = x^3y + xy^4}$.}
\[\J = (3x^2y + y^4,\,x^3+4xy^3)
\phantom{XXX} q_x =\tfrac{3}{11},\,q_y=\tfrac{2}{11},\,\hat c=\tfrac{12}{11}
\phantom{XXX} G = \br{J} = \Z/11\Z\]
\[\Fix J^k = \begin{cases}
  \C^2_{xy} & \text{if } k=0\\
  0 & \text{otherwise}
\end{cases}
\phantom{XXX}
\Q|_{\Fix J^k} = \begin{cases}
  \br{1,\,y,\,x,\,y^2,\,xy,\,y^3,\,x^2,\,xy^2,\,x^2y,\,xy^3,\,x^2y^2,\,x^2y^3},\, \mu=12\\
  \br{1}
\end{cases} \]

\begin{center}
\begin{tabular}{|c||c|c|c|c|c|c|c|c|c|c|c|}\hline
$k$	&0	&1	&2	&3	&4	&5	&6	&7	&8	&9	&10\\\hline
$|G|\cdot\deg_W$	&12	&0	&10	&20	&8	&18	&6	&16	&4	&14	&24\\\hline
invariants	& $x^2e_0,\,y^3e_0$	& $\1$	& $e_{2}$	& $e_{3}$	& $e_{4}$	& $e_{5}$	& $e_{6}$	& $e_{7}$	& $e_{8}$	& $e_{9}$	& $e_{10}$\\\hline
\end{tabular}
\end{center}

Potential non-zero correlators:

By the pairing axiom,
$\br{\1,	\1, e_{10}}$,
$\br{\1,	e_2, e_{9}}$,
$\br{\1,	e_3, e_{8}}$,
$\br{\1,	e_{4},	e_{7}}$, and
$\br{\1,	e_{5},	e_{6}}$,
$\br{\1,	e_{9},	e_{9}}$ are all equal to $1$.

Also, 
$\br{\1,x^2e_0,x^2e_0}=\tfrac{1}{11}$, $\br{\1,x^2e_0,y^3e_0}=-\tfrac{4}{11}$ and $\br{\1,y^3e_0,y^3e_0}=-\tfrac{3}{11}$.


By the concavity axiom,
$\br{e_{2},	e_{2},	e_{8}}$,
$\br{e_{2},	e_{4},	e_{6}}$,
$\br{e_{6},	e_{8},	e_{9}}$, and
$\br{e_{7},	e_{8},	e_{8}}$
are all equal to $1$.

By the index-zero axiom,
$\br{e_{4}, e_{4}, e_{4}} = {-3}$.

The correlators $\br{x^2e_0,e_{4},e_{8}}$, $\br{y^3e_0,e_{4},e_{8}}$, $\br{x^2,e_{6},e_{6}}$ and $\br{y^3e_0,e_{6},e_{6}}$ are also non-zero, and we use the composition axiom to extract the necessary information (although we avoid computing individual correlators directly).

Note
\begin{gather*}
e_6^2 = \sum_{\mu,\,\nu\in\{x^2 e_0,\,y^3 e_0\}} \br{e_{6},	e_{6}, \mu}\eta^{\mu,\nu}\nu\\
e_6^3 = \sum_{\mu,\,\nu\in\{x^2 e_0,\,y^3 e_0\}} \br{e_{6},	e_{6}, \mu}\eta^{\mu,\nu} \br{\nu, e_6, e_6} e_5,\\
\end{gather*}
By the composition axiom, the coefficient of $e_5$ here is just the value of the four pointed class
$\Lambda_{0,4}^{Z_{12}}(e_6,e_6,e_6,e_6)$.
This four-pointed class has codimension zero, and $\lb_x = -2$ and $\lb_y = 0$,
so its value is the $y$-degree of the Witten map, which is ${-4}$.

Consider the map
\[\varphi:\C[X,Y]\to \FJRW{Z_{12}}{G}\]
defined by $X\mapsto e_{6}$, and $Y\mapsto e_{8}$ extending by $\C$-linearity and multiplicativity.  One can check directly that this maps onto each twisted sector.  

To prove surjectivity, we need to check that it also maps onto the untwisted sector.  This is equivalent to the linear independence of $e_6^2$ and $e_8^3$.  We note the following identities, which are direct consequences of the three-point correlator values cited above:
\begin{gather*}
  e_6^2*e_6 = -4e_5\\
  e_8^3*e_6 = e_8^2*(e_8e_6) = e_4*e_2=e_5\\
  \\
  e_6^2*e_8= e_6*(e_6e_8) = e_6*e_2 = e_7\\
  e_8^3*e_8 = e_8^2*e_8^2 = e_4*e_4 = -3e_7.
\end{gather*}

Putting $\mu = e_6^2+4e_8^3$ and $\nu = 3e_6 + e_8^3$, we see that
\begin{gather*}
\mu*e_6 =0\\
\mu*e_8 = -11e_7\\
\\
\nu*e_6 = -11e_5\\
\nu*e_8 = 0,
\end{gather*}
and conclude that $\mu$ and $\nu$ are linearly independent combinations of $e_6^2$ and $e_8^3$, yielding surjectivity of the map $\varphi$.

The above identities also show $(3X^2Y+Y^4,\,X^3+4XY^3)\subseteq \ker \varphi$. Since
$\C[X,Y]/(3X^2Y+Y^4,\,X^3+4XY^3)=\q{Z_{12}^T}$ has dimension 12, we deduce the inclusion is in fact equality, and $\FJRW{Z_{12}}{G}\cong\q{Z_{12}^T}$.
\newline\hrule

\subsection*{$\mathbf{Z_{13} = x^3y + y^6}$.}
\[\J = (3x^2y,\,x^3+6y^5)
\phantom{XXX} q_x =\tfrac{5}{18},\,q_y=\tfrac{3}{18},\,\hat c=\tfrac{20}{18}
\phantom{XXX} G = \br{J} = \Z/18\Z\]
\[\Fix J^k = \begin{cases}
  \C^2_{xy} & \text{if } k=0\\
  \C_{y} & \text{if } k=6,\,12\\
  0 & \text{otherwise}
\end{cases}
\phantom{XXX}
\Q|_{\Fix J^k} = \begin{cases}
  \br{1,\,y,\,x,\,y^2,\,xy,\,y^3,\,x^2,\,xy^2,\,y^4,\,xy^3,\,y^5,\,xy^4,\,xy^5},\, \mu=13\\
  \br{1,\,y,\,y^2,\,y^3,\,y^4},\,\mu=5\\
  \br{1}\\
\end{cases} \]

\begin{center}
\begin{tabular}{|c||c|c|c|c|c|c|c|c|c|c|c|c|c|c|c|c|}\hline
$k$	&0	&1	&2	&3	&4	&5	&7	&8	&9	&10	&11	&13	&14	&15	&16	&17\\\hline
$|G|\cdot\deg_W$	&20	&0	&16	&32	&12	&28	&24	&4	&20	&36	&16	&12	&28	&8	&24	&40\\\hline
invariants	& $x^2e_{0}$	&$\1$	&$e_{2}$	&$e_{3}$	&$e_{4}$	&$e_{5}$	&$e_{7}$	&$e_{8}$	&$e_{9}$	&$e_{10}$	&$e_{11}$	 &$e_{13}$	&$e_{14}$	&$e_{15}$	&$e_{16}$	&$e_{17}$\\\hline
\end{tabular}
\end{center}

Potential non-zero correlators:

By the pairing axiom,
$\br{\1,	\1, e_{17}}$,
$\br{\1,	e_2, e_{16}}$,
$\br{\1,	e_3, e_{15}}$,
$\br{\1,	e_{4},	e_{14}}$,
$\br{\1,	e_{5},	e_{13}}$,
$\br{\1,	e_{7},	e_{11}}$,
$\br{\1,	e_{8},	e_{10}}$,
$\br{\1,	e_{9},	e_{9}}$ are all equal to $1$, and
$\br{\1, x^2e_{0}, x^2e_{0}}={-\tfrac{1}{3}}$.

By the concavity axiom,
$\br{e_{2},	e_{2},	e_{15}}$,
$\br{e_{2},	e_{4},	e_{13}}$,
$\br{e_{2},	e_{8},	e_{9}}$,
$\br{e_{3},	e_{8},	e_{8}}$,
$\br{e_{4},	e_{7},	e_{8}}$,
$\br{e_{7},	e_{15},	e_{15}}$,
$\br{e_{8},	e_{13},	e_{16}}$,
$\br{e_{8},	e_{14},	e_{15}}$,
$\br{e_{9},	e_{13},	e_{15}}$, and
$\br{e_{11},	e_{13},	e_{13}}$
are all equal to $1$.

By the index-zero axiom,
$\br{e_{4}, e_{4}, e_{11}}$ and
$\br{e_{11}, e_{11}, e_{15}}$ are both equal to ${-3}$.

The composition axiom can be used to compute
$\br{x^2e_{0},	e_{4},	e_{15}}$, and
$\br{x^2e_{0},	e_{8},	e_{11}}$, but we do not need these to establish the desired isomorphism.

Consider the map
\[\varphi:\C[X,Y]\to \FJRW{Z_{13}}{G},\]
defined by $X\mapsto e_{13}$, and $Y\mapsto e_{8}$ extending by $\C$-linearity and multiplicativity.  One can check directly that this map is surjective.

Our correlators tell us that $e_{8}^3*e_{8}^3 = -3e_{7} = -3e_{13}^2$ and $e_{13}e_{8}^5 = 0$.

Hence $(3X^2+Y^6,\,6XY^5)\subseteq \ker \varphi$.  Since
$\C[X,Y]/(3X^2+Y^6,\,6XY^5)=\q{Z_{13}^T}$ has dimension 16, we deduce the inclusion is in fact equality, and $\FJRW{Z_{13}}{G}\cong\q{Z_{13}^T}$.
\newline\hrule

\subsection*{$\mathbf{W_{12} = x^4 + y^5}$.}
By Axiom \ref{ax:sums}, the conjecture holds for $W_{12}$ because it holds for the simple singularities $A_3 = x^4$ and $A_4 = y^5$.\newline\hrule

\subsection*{$\mathbf{W_{13} = x^4 + xy^4}$.}
\[\J = (4x^3+y^4,\,4xy^3)
\phantom{XXX} q_x =\tfrac{4}{16},\,q_y=\tfrac{3}{16},\,\hat c=\tfrac{18}{16}
\phantom{XXX} G = \br{J} = \Z/16\Z\]
\[\Fix J^k = \begin{cases}
  \C^2_{xy} & \text{if } k=0\\
  \C_{x} & \text{if } k=4,\,8,\,12\\
  0 & \text{otherwise}
\end{cases}
\phantom{XXX}
\Q|_{\Fix J^k} = \begin{cases} \br{1,\,y,\,x,\,y^2,\,xy,\,x^2,\,y^3,\,xy^2,\,x^2y,\,y^4,\,y^5,\,y^6},\, \mu=13\\
  \br{1,\,x,\,x^2},\,\mu=3\\
  \br{1}\\
\end{cases}\]

\begin{center}
\begin{tabular}{|c||c|c|c|c|c|c|c|c|c|c|c|c|c|}\hline
$k$	&0	&1	&2	&3	&5	&6	&7	&9	&10	&11	&13	&14	&15\\\hline
$|G|\cdot\deg_W$	&18	&0	&14	&28	&24	&6	&20	&16	&30	&12	&8	&22	&36\\\hline
invariants	& $y^3e_{0}$	& $ \1 $	& $e_{2}$	& $e_{3}$	& $e_{5}$	& $e_{6}$	& $e_{7}$	& $e_{9}$	& $e_{10}$	& $e_{11}$	& $e_{13}$	& $e_{14}$	& $e_{15}$\\\hline
\end{tabular}
\end{center}

Potential non-zero correlators:

By the pairing axiom,
$\br{\1, \1, e_{15}}$,
$\br{\1,	e_2, e_{14}}$,
$\br{\1,	e_3, e_{13}}$,
$\br{\1,	e_5, e_{11}}$,
$\br{\1,	e_{6},	e_{10}}$, and
$\br{\1,	e_{7},	e_{9}}$ are all equal to $1$, and
$\br{\1,	y^3e_0, y^3e_{0}}={-\tfrac{1}{4}}$

By the concavity axiom,
$\br{e_{2},	e_{2},	e_{13}}$,
$\br{e_{2},	e_{6},	e_{9}}$,
$\br{e_{5},	e_{6},	e_{6}}$,
$\br{e_{6},	e_{13},	e_{14}}$,
$\br{e_{7},	e_{13},	e_{13}}$, and
$\br{e_{9},	e_{11},	e_{13}}$
are all equal to $1$.

By the index-zero axiom, $\br{e_{11},	e_{11}, e_{11}}={-4}$.

The correlator $\br{y^3e_0,	e_6, e_{11}}$ may be non-zero, and the composition axiom can be used to compute it.  However, using associativity of the product we can compute the ring structure on $\FJRW{W_{13}}{G}$ without it.

Consider the map
\[\varphi:\C[X,Y]\to \FJRW{W_{13}}{G}\]
defined by $X\mapsto e_{6}$, and $Y\mapsto e_{13}$, extending by $\C$-linearity and multiplicativity.  One can check directly that this map is surjective.

Our correlators tell us that $e_{6}^2*e_6^2 = -4e_5 = -4e_{13}^3$, and $e_6^2*(e_6e_{13}) = 0$.

Hence $(4X^3Y,\, X^4+4Y^3)\subseteq \ker \varphi$.  Since
$\C[X,Y]/(4X^3Y,\, X^4+4Y^3)=\q{W_{13}^T}$ has dimension 13, we deduce the inclusion is in fact equality, and $\FJRW{W_{13}}{G}\cong\q{W_{13}^T}$.
\newline\hrule

\subsection*{$\mathbf{E_{12} = x^3 + y^7}$.}

By Axiom \ref{ax:sums}, the Mirror Symmetry Conjecture holds for $E_{12}$ because it holds for the simple singularities $A_2=x^3$ and $A_6=y^7$.\newline\hrule

\subsection*{$\mathbf{E_{13} = x^3 + xy^5}$.}
\[\J = (4x^2+y^5,\,5xy^4)
\phantom{XXX} q_x =\tfrac{5}{15},\,q_y=\tfrac{2}{15},\,\hat c=\tfrac{16}{15}
\phantom{XXX} G = \br{J} = \Z/15\Z\]
\[\Fix J^k = \begin{cases}
  \C^2_{xy} & \text{if } k=0\\
  \C_{x} & \text{if } k=3,\,6,\,9,\,12\\
  0 & \text{otherwise}
\end{cases}
\phantom{XXX}
\Q|_{\Fix J^k} = \begin{cases}
  \br{1,\,y,\,y^2,\,x,\,y^3,\,xy,\,y^4,\,xy^2,\,y^5,\,xy^3,\,y^6,\,y^7,\,y^8},\, \mu=13\\
  \br{1,\,x},\,\mu=2\\
  \br{1}\\
\end{cases}\]

\begin{center}
\begin{tabular}{|c||c|c|c|c|c|c|c|c|c|c|c|}\hline
$k$	&0	&1	&2	&4	&5	&7	&8	&10	&11	&13	&14\\\hline
$|G|\cdot\deg_W$	&16	&0	&14	&12	&26	&24	&8	&6	&20	&18	&32\\\hline
invariants	& $y^4e_{0}$	& $\1$ 	& $e_{2}$ 	& $e_{4}$ 	& $e_{5}$ 	& $e_{7}$ 	& $e_{8}$ 	& $e_{10}$ 	& $e_{11}$ 	& $e_{13}$ 	& $e_{14}$ \\\hline
\end{tabular}
\end{center}

Potential non-zero correlators:

By the pairing axiom,
$\br{\1, \1, e_{14}}$,
$\br{\1,	e_2, e_{13}}$,
$\br{\1,	e_4, e_{11}}$,
$\br{\1,	e_5, e_{10}}$, and
$\br{\1,	e_{7},	e_{8}}$ are all equal to $1$, and
$\br{\1,	y^4e_{0},	y^4e_{0}}=-\tfrac{1}{5}$.

By the concavity axiom,
$\br{e_{2},	e_{4},	e_{10}}$,
$\br{e_{4},	e_{4},	e_{8}}$,
$\br{e_{8},	e_{10},	e_{13}}$, and
$\br{e_{10},	e_{10},	e_{11}}$
are all equal to $1$.

By the composition axiom, we have
\[-5 = \br{e_{8},	e_{8}, y^4e_0}\eta^{y^4e_0,y^4e_0}\br{y^4e_0,e_8,e_8},\]
so $\br{e_{8},	e_{8}, y^4e_0}=\pm 1$.

Consider the map
\[\varphi:\C[X,Y]\to \FJRW{E_{13}}{G}\]
defined by $X\mapsto e_{8}$, and $Y\mapsto e_{10}$, extending by $\C$-linearity and multiplicativity.  One can check directly that this map is surjective.

Our correlators tell us that $e_{8}*e_8e_{10} = 0$ and $e_{8}^2*e_{8} = -5e_{7} = -5e_{10}^4$.

Hence $(3X^2Y,\, X^3+5Y^4)\subseteq \ker \varphi$. Since
$\C[X,Y]/(3X^2Y,\, X^3+5Y^4)=\q{E_{13}^T}$ has dimension 11, we deduce the inclusion is in fact equality, and $\FJRW{E_{13}}{G}\cong\q{E_{13}^T}$.
\newline\hrule

\subsection*{$\mathbf{E_{14} = x^3 + y^8}$.}
By Axiom \ref{ax:sums}, the Mirror Symmetry Conjecture holds for $E_{14}$ because it holds for the simple singularities $A_2=x^3$ and $A_7=y^8$.

\subsection{Bimodal Singularities}\label{bi}

We now turn to the fourteen exceptional bimodal families listed in Arnol'd \cite{Arn-1}. We note that the conjecture has already been shown for $E_{18}$, $E_{20}$, $U_{16}$, $W_{18}$, $Q_{16}$, and $Q_{18}$ in \cite{FJR} as these are sums of simple singularities. Also $E_{19}$ was constructed in the introduction, so in this section, we will consider the eight remaining families, together with their mirror partners.

In \cite{FJR}, the construction of the FJRW-ring does not include singularities with a weight greater than or equal to 1/2. In particular, it has not been verified that the required compact moduli space exists in these cases, and that the three-point correlators satisfy the proper axioms. However, Fan-Jarvis-Ruan have proven that the construction works for two particular singularities having a variable with weight 1/2, namely $A_2$ and $D_n^T$ for $n$ even, and it is expected to work for all cases with weight $\frac12$ (see \cite{fjr2}). When our singularities have variables of weight $\frac12$, we will treat them assuming the FJR axioms hold.

In this paper, we have included four examples of singularities with this property, namely $Q_{17}^T$, $S_{17}^T$, $Q_{2,0}^T$ and $S_{1,0}^T$. In each case, we have shown how the isomorphism between the FJRW-rings and the Milnor rings should work under the assumption that a certain three-point correlator is non-zero.

\subsection*{$\mathbf{Z_{17}=x^3y +y^8}$}
\[
\jac=(3x^2y,x^3+8y^7)
\phantom{XXX}
q_x = \tfrac 7{24},\ \ q_y = \tfrac 3{24},\ \ \hat c=\tfrac {28}{24}
\phantom{XXX}
G=\br J\cong \Z/24\Z
\]

\[
\Fix {J^k} = \begin{cases}
    \C^2 & \text{if}\quad k = 0\\
    \C_y & \text{if}\quad 8|k, k\neq 0\\
   0    & \text{otherwise}
\end{cases}
\phantom{XXX}
\q{}|_{\Fix J^k} =
\begin{cases}
\br{1,x,x^2, y, \dots, y^7, xy, \dots,xy^7},\ \mu=17\\
\br{1, y, \dots, y^6},\ \mu=7\\
\br{1} 
\end{cases}
\]

\begin{center}
\begin{tabular}{|c||c|c|c|c|c|c|c|c|c|c|c|}
\hline
k       & 0&1&2&3&4&5&6&7&9&10&11\\
\hline
$|G|\cdot\deg_W$  & 28  &0 &20  &40 &12 &32 &52 &24   &16 &36 &8\\
\hline
invariants &$x^2e_0$& $\1$& $e_2$& $e_3$& $e_4$& $e_5$& $e_6$& $e_7$&  $e_9$& $e_{10}$& $e_{11}$\\
\hline
\end{tabular}

\vspace*{\baselineskip}

\begin{tabular}{|c||c|c|c|c|c|c|c|c|c|c|c|}
\hline
k       & 12 & 13 & 14& 15  & 17 & 18 & 19  & 20 & 21   & 22 & 23\\
\hline
$|G|\cdot\deg_W$  &28 &48 &20 &12 &32 &4 &24 &44 &16 &36 &56\\
\hline
invariants & $e_{12}$& $e_{13}$&  $e_{14}$& $e_{15}$& $e_{17}$& $e_{18}$& $e_{19}$&  $e_{20}$& $e_{21}$&  $e_{22}$& $e_{23}$\\
\hline
\end{tabular}
\end{center}

Potential non-zero correlators:

By the pairing axiom $\br{\1,\1,e_{23}}$, $\br{\1,e_2,e_{22}}$, $\br{\1,e_3,e_{21}}$, $\br{\1,e_4,e_{20}}$, $\br{\1,e_5,e_{19}}$, $\br{\1,e_6,e_{18}}$, $\br{\1,e_7,e_{17}}$, $\br{e_1,e_9,e_{15}}$, $\br{\1,e_{10},e_{14}}$, $\br{\1,e_{11},e_{13}}$ and $\br{\1,e_{12},e_{12}}$ are equal to 1, and $\br{x^2e_0,x^2e_0,\1}=-1/3$.

By the Concavity axiom $\br{e_2,e_2,e_{21}}$, $\br{e_2,e_4,e_{19}}$, $\br{e_2,e_5,e_{18}}$, $\br{e_2,e_9,e_{14}}$, $\br{e_2,e_{11},e_{12}}$, $\br{e_3,e_4,e_{18}}$, $\br{e_3,e_{11},e_{11}}$, $\br{e_4,e_4,e_{17}}$, $\br{e_4,e_9,e_{12}}$, $\br{e_4,e_{10},e_{11}}$, $\br{e_5,e_9,e_{11}}$, $\br{e_7,e_9,e_9}$, $\br{e_9,e_{18},e_{22}}$, $\br{e_9,e_{19},e_{21}}$, $\br{e_{10},e_{18},e_{21}}$, $\br{e_{11},e_{17},e_{21}}$, $\br{e_{11},e_{18},e_{20}}$, $\br{e_{11},e_{19},e_{19}}$, $\br{e_{12},e_{18},e_{19}}$, $\br{e_{13},e_{18},e_{18}}$, and $\br{e_{14},e_{17},e_{18}}$ are all equal to 1.

By the Index Zero axiom $\br{e_4,e_7,e_{14}}=\br{e_7,e_7,e_{11}}=\br{e_7,e_{21},e_{21}}= \br{e_{14},e_{14},e_{21}}=-3$.

By the composition axiom $\br{x^2e_0,e_4,e_{21}}=\br{x^2e_0,e_7,e_{18}}=\br{x^2e_0,e_{11},e_{14}}=\pm 1$.

Consider the map $\varphi:\C[X,Y]\rightarrow \h{Z_{17}}{G}$ defined by $X\mapsto e_9$ and $Y\mapsto e_{18}$ and extending to a $\C$-algebra homomorphism. One can check directly that this map is surjective. 

Our correlators tell us that $e_9^2=e_{17}$, $e_{18}^8=-3e_{17}$, and $e_9e_{18}^7=0$. 

Hence $(3X^2+Y^8, 8XY^8)\subset\ker\varphi$. Since $\C[X,Y]/(3X^2+Y^8, 8XY^8)=\q{Z_{17}^T}$ has dimension 22, we deduce that the inclusion is equality, and so we have the isomorphism $\h{Z_{17}}{G}\cong\q{Z_{17}^T}$. 
\newline\hrule

\subsubsection*{$\mathbf{Z_{17}^T=x^3 + xy^8}$}
\[
\jac=(3x^2+y^8, 8xy^7)
\phantom{XXX}
q_x = \tfrac 8{24}, q_y = \tfrac 2{24},\hat c=\tfrac {28}{24}
\phantom{XXX}
G\cong \Z/24\Z
\phantom{XXX}
\br{J}\cong \Z/12\Z
\]
Since $J$ does not generate the maximal group of symmetries, we will use the generator $\g=(\zeta^{16},\zeta)$, where $\zeta^{24}=1$. We will index our graded FJRW-ring by powers of $\g$.
\[
\Fix {\g^k} = \begin{cases}
    \C^2 & \text{if}\quad k = 0\\
    \C_x & \text{if}\quad 3|k, k\neq 0\\
   0    & \text{otherwise}
\end{cases}
\phantom{XXX}
\q{}|_{\Fix \g^k} =
\begin{cases}
\br{1,x,x^2, y,\dots, y^{7}, xy,\dots, xy^6,x^2y,\dots,x^2y^6},\,\mu=22\\
\br{1,x},\,\mu=2\\
\br{1} 
\end{cases}
\]

\begin{center}
\begin{tabular}{|c||c|c|c|c|c|c|c|c|c|c|c|c|c|c|c|c|c|}
\hline
k       & 0&1&2&4&5&7&8&10&11& 13 & 14& 16  & 17 & 19  & 20 & 22 & 23\\
\hline
$|G|\cdot\deg_W$  &28   &14 &0 &20  &6 &26 &12 &32 &18 &38 &24  &44  &30 &50 &36 &56 &42\\
\hline
invariants & $y^7e_0$& $e_1$& $\1$& $e_4$& $e_5$& $e_7$& $e_8$& $e_{10}$& $e_{11}$& $e_{13}$& $e_{14}$& $e_{16}$& $e_{17}$& $e_{19}$& $e_{20}$& $e_{22}$& $e_{23}$\\
\hline
\end{tabular}
\end{center}

Potential non-zero correlators:

By the pairing axiom $\br{\1,e_1,e_{23}}$, $\br{\1,\1,e_{22}}$, $\br{\1,e_4,e_{20}}$, $\br{\1,e_5,e_{19}}$, $\br{\1,e_7,e_{17}}$, $\br{\1,e_8,e_{16}}$, $\br{\1,e_{10},e_{14}}$, and $\br{\1,e_{11},e_{13}}$ are equal to 1, and $\br{y^7e_0,y^7e_0,\1}=-1/8$.

By the Concavity axiom $\br{e_1,e_5,e_{20}}$, $\br{e_1,e_8,e_{17}}$, $\br{e_1,e_{11},e_{14}}$, $\br{e_4,e_5,e_{17}}$, $\br{e_4,e_{8},e_{14}}$, $\br{e_4,e_{11},e_{11}}$, $\br{e_5,e_{5},e_{16}}$, $\br{e_5,e_7,e_{14}}$, $\br{e_5,e_8,e_{13}}$, $\br{e_5,e_{10},e_{11}}$, $\br{e_7,e_8,e_{11}}$, and $\br{e_8,e_8,e_{10}}$ are all equal to 1.

By the composition axiom $\br{e_1,e_1,y^7e_{0}}=\pm 1$.

Consider the map $\varphi:\C[X,Y]\rightarrow \h{Z_{17}^T}{G}$ defined by $X\mapsto e_1$ and $Y\mapsto e_5$ and extending to a $\C$-algebra homomorphism. One can check directly that this map is surjective. 

Our correlators tell us that $e_1^3=-8e_{23}$, $e_5^7=e_{23}$, and $e_1^2e_5=0$. 

Hence $(3X^2Y, X^3+8Y^7)\subset\ker\varphi$. Since $\C[X,Y]/(3X^2Y, X^3+8Y^7)=\q{Z_{17}}$ has dimension 17, we deduce that the inclusion is equality, and so we have the isomorphism $\h{Z_{17}^T}{G}\cong\q{Z_{17}}$. 
\newline\hrule

\subsubsection*{$\mathbf{Z_{18}=x^3y +xy^6}$}
\[
\jac=(3x^2y+y^6, x^3+6xy^5)
\phantom{XXX}
q_x=\tfrac 5{17}, q_y=\tfrac 2{17}, \hat c=\tfrac {20}{17}
\phantom{XXX}
G=\br J\cong \Z/17\Z
\]
\[
\Fix {J^k} = \begin{cases}
        \C^2 & \text {if}\quad k = 0\\
        0    & \text {otherwise}
\end{cases} 
\phantom{XXX}
\q{}|_{\Fix J^0} =
\begin{cases}
\br{1,x,x^2,y,\dots, y^{7},xy, \dots, xy^5,x^2y,\dots,x^2y^6},\mu=18\\
\br{1} 
\end{cases}
\]

\begin{center}
\begin{tabular}{|c||c|c|c|c|c|c|c|c|c|c|c|c|c|c|c|c|c|}
\hline
k       & 0    &1 &2     &3     &4    &5     &6     &7     &8 & 9    &10    &11   &12    &13    &14    &15    &16\\
\hline
$|G|\cdot\deg_W$  &20 &0 &14 &28 &8 &22 &36 &16 &30 &10 &24 &4 &18 &32 &12 &26 &40\\
\hline
invariants &$x^2e_0, y^5e_0$& $\1$& $e_2$& $e_3$& $e_4$& $e_5$& $e_6$& $e_7$& $e_8$& $e_9$& $e_{10}$& $e_{11}$& $e_{12}$& $e_{13}$& $e_{14}$& $e_{15}$& $e_{16}$\\
\hline
\end{tabular}
\end{center}

Potential non-zero correlators:

By the pairing axiom $\br{\1,\1,e_{16}}$, $\br{\1,e_2,e_{15}}$, $\br{\1,e_3,e_{14}}$, $\br{\1,e_4,e_{13}}$, $\br{\1,e_5,e_{12}}$, $\br{\1,e_6,e_{11}}$, $\br{\1,e_7,e_{10}}$, and $\br{\1,e_8,e_9}$ are equal to 1, and $\br{x^2e_0,x^2e_0,\1}=-6/17$, $\br{x^2e_0,y^5e_0,\1}=1/17$, $\br{y^5e_0,y^5e_0,\1}=-3/17$.

By the Concavity axiom $\br{e_2,e_2,e_{14}}$, $\br{e_2,e_4,e_{12}}$, $\br{e_2,e_5,e_{11}}$, $\br{e_2,e_7,e_9}$, $\br{e_3,e_4,e_{11}}$, $\br{e_4,e_4,e_{10}}$, $\br{e_4,e_5,e_9}$, $\br{e_9,e_{11},e_{15}}$, $\br{e_9,e_{12},e_{14}}$, $\br{e_{10},e_{11},e_{14}}$, $\br{e_{11},e_{11},e_{13}}$, and $\br{e_{11},e_{12},e_{12}}$ are all equal to 1.

By the Index Zero axiom $\br{e_{4},e_{7},e_{7}}=\br{e_{7},e_{14},e_{14}}=-3$

The remaining potentially non-zero three point correlators cannot be determined from the axioms alone. These correlators are $\br{x^2e_{0},e_{4},e_{14}}$, $\br{y^5e_{0},e_{4},e_{14}}$, $\br{x^2e_{0},e_{7},e_{11}}$, $\br{y^5e_{0},e_{7},e_{11}}$, $\br{x^2e_{0},e_{9},e_{9}}$, and $\br{y^5e_{0},e_{9},e_{9}}$. 

In order to prove this we consider the homomorphim $\varphi:\C[X,Y]\rightarrow \h{Z_{18}}{G}$ defined by $X\mapsto e_9$ and $Y\mapsto e_{11}$. It takes some work to show that $\varphi$ is surjective in this case. One can check that $\varphi$ maps \emph{onto} the one-dimensional sectors. 

However, the sector corresponding to the identity in $G$ is two-dimensional. To prove surjectivity, it suffices to show that $e_{11}^5$ and $e_9^2$ are linearly independent.

We note the following identities, which are direct consequences of the three-point correlator values cited above:
\begin{gather*}
  e_9^2*e_9 = -6e_8\\
  e_9^2*e_{11} = e_{10}\\
  \\
  e_{11}^5*e_9= e_8\\
  e_{11}^5*e_{11} = -3e_{10}.
\end{gather*}
Putting $\mu = e_9^2+6e_{11}^5$ and $\nu = 3e_9^2 + e_{11}^5$, we see that
\begin{gather*}
\mu*e_9 =0\\
\mu*e_{11} = -17e_{10}\\
\\
\nu*e_9 = -17e_8\\
\nu*e_{11} = 0,
\end{gather*}
and conclude that $\mu$ and $\nu$ are linearly independent combinations of $e_9^2$ and $e_{11}^5$, yielding surjectivity of the map $\varphi$.

The above identities also tell us that  $(3X^2Y+Y^6, X^3+6XY^5)\subset \ker\varphi$. Since $\C[X,Y]/(3X^2Y+Y^6, X^3+6XY^5)=\q{Z_{18}}$ has dimension 18 we deduce that we have an isomorphism $\q{Z_{18}}\cong\h{Z_{18}}{G}$.
\newline\hrule

\subsubsection*{$\mathbf{Z_{19}=x^3y +y^9}$}
\[
\jac=(3x^2y, x^3+9y^8)
\phantom{XXX}
q_x=\tfrac 8{27}, q_y=\tfrac 3{27}, \hat c=\tfrac {32}{27}
\phantom{XXX}
G=\br J\cong \Z/27\Z
\]
\[
\Fix {J^k} = \begin{cases}
   \C^2_{xy} & \text {if}\quad k = 0\\
    \C_y & \text {if}\quad 9|k, k\neq 0\\
   0    & \text {otherwise}
\end{cases}
\phantom{XXX}
\q{}|_{\Fix J^k} =
\begin{cases}
\br{1,x,x^2,y, \dots, y^8,xy, xy^2, \dots, xy^8},\mu=19\\
\br{1,y,y^2,y^3,y^4,y^5,y^6,y^7}\mu=8\\
\br{1} 
\end{cases}
\]

\begin{center}
\begin{tabular}{|c||c|c|c|c|c|c|c|c|c|c|c|c|c|c|c|c|c|}
\hline
k       & 0     &1 &2     &3     &4    &5     &6     &7     &8 & 10    &11   &12    &13    &14    &15    &16 & 17\\
\hline
$|G|\cdot\deg_W$  &32 &0 &22 &44 &12 &34 &56 &24 &46 &36 &4 &26 &48 &16 &38 &60 &28\\
\hline
invariants & $x^2e_0$& $\1$& $e_2$& $e_3$& $e_4$& $e_5$& $e_6$& $e_7$& $e_8$& $e_{10}$& $e_{11}$& $e_{12}$& $e_{13}$& $e_{14}$& $e_{15}$& $e_{16}$& $e_{17}$\\
\hline
\end{tabular}

\vspace*{\baselineskip}

\begin{tabular}{|c||c|c|c|c|c|c|c|c|}
\hline
k       & 19 & 20   &21   & 22  & 23   &24    &25   &26 \\
\hline
$|G|\cdot\deg_W$  &18 &40 &8 &30 &52 &20 &42 &64\\
\hline
invariants & $e_{19}$& $e_{20}$& $e_{21}$& $e_{22}$& $e_{23}$& $e_{24}$& $e_{25}$& $e_{26}$\\
\hline
\end{tabular}
\end{center}

Potential non-zero correlators:

By the pairing axiom $\br{\1,\1,e_{26}}$, $\br{\1,e_2,e_{25}}$, $\br{\1,e_3,e_{24}}$, $\br{\1,e_4,e_{23}}$, $\br{\1,e_5,e_{22}}$, $\br{\1,e_6,e_{21}}$, $\br{\1,e_7,e_{20}}$, $\br{\1,e_8,e_{19}}$, $\br{\1,e_{10},e_{17}}$, $\br{\1,e_{11},e_{16}}$, $\br{\1,e_{12},e_{15}}$,and $\br{\1,e_{13},e_{14}}$ are equal to 1, and $\br{x^2e_0,x^2e_0,\1}=-1/3$,

By the Concavity axiom $\br{e_2,e_2,e_{24}}$, $\br{e_2,e_4,e_{22}}$, $\br{e_2,e_5,e_{21}}$, $\br{e_2,e_7,e_{19}}$, $\br{e_2,e_{11},e_{15}}$, $\br{e_2,e_{12},e_{14}}$, $\br{e_3,e_4,e_{21}}$, $\br{e_3,e_{11},e_{14}}$, $\br{e_4,e_4,e_{20}}$, $\br{e_4,e_5,e_{19}}$, $\br{e_4,e_{10},e_{14}}$, $\br{e_4,e_{11},e_{13}}$, $\br{e_4,e_{12},e_{12}}$, $\br{e_5,e_{11},e_{12}}$, $\br{e_6,e_{11},e_{11}}$, $\br{e_7,e_{10},e_{11}}$, $\br{e_{10},e_{21},e_{24}}$, $\br{e_{11},e_{19},e_{25}}$, $\br{e_{11},e_{20},e_{24}}$, $\br{e_{11},e_{21},e_{23}}$, $\br{e_{11},e_{22},e_{22}}$, $\br{e_{12},e_{19},e_{24}}$, $\br{e_{12},e_{21},e_{22}}$, $\br{e_{13},e_{21},e_{21}}$, $\br{e_{14},e_{17},e_{24}}$, $\br{e_{14},e_{19},e_{22}}$, $\br{e_{14},e_{20},e_{21}}$,  $\br{e_{15},e_{19},e_{21}}$, and $\br{e_{17},e_{19},e_{19}}$ are all equal to 1.

By the Index Zero axiom $\br{e_{4},e_{7},e_{17}}=\br{e_{7},e_{7},e_{14}}=\br{e_{7},e_{24},e_{24}}=\br{e_{17},e_{17},e_{21}}=-3$

By the composition axiom $\br{x^2e_{0},e_{4},e_{24}}=\br{x^2e_{0},e_{7},e_{21}}=\br{x^2e_{0},e_{11},e_{17}}=\br{x^2e_{0},e_{14},e_{14}}=\pm 1$.

Consider the map $\varphi:\C[X,Y]\rightarrow \h{Z_{19}}{G}$ defined by $X\mapsto e_{19}$ and $Y\mapsto e_{11}$ and extending to a $\C$-algebra homomorphism. One can check directly that this map is surjective. 

Our correlators tell us that $e_{19}^2=e_{10}$, $e_{11}^9=-3e_{10}$ and $e_{19}e_{11}^8=0$. 

Hence $(3X^2+Y^9, 9XY^8)\subset\ker\varphi$. Since $\C[X,Y]/(3X^2+Y^9, 9XY^8)=\q{Z_{19}^T}$ has dimension 25, we deduce that the inclusion is equality, and so we have the isomorphism $\h{Z_{19}}{G}\cong\q{Z_{19}^T}$.
\newline\hrule

\subsubsection*{$\mathbf{Z_{19}^T=x^3 + xy^9}$}
\[
\jac=(3x^2+y^9, 9xy^8)
\phantom{XXX}
q_x=\tfrac 9{27}, q_y=\tfrac 2{27}, \hat c=\tfrac {32}{27}
\phantom{XXX}
G=\br J\cong \Z/27\Z
\]
\[
\Fix {J^k} = \begin{cases}
   \C^2 & \text {if}\quad k = 0\\
    \C_x & \text {if}\quad 3|k\\
   0    & \text {otherwise}
\end{cases}
\phantom{XXX}
\q{}|_{\Fix J^k} =
\begin{cases}
\br{1,x,y,\dots, y^{8},xy, \dots, xy^7,x^2y,\dots,x^2y^7},\mu=25\\
\br{1,x},\mu=2\\
\br{1} 
\end{cases}
\]

\begin{center}
\begin{tabular}{|c||c|c|c|c|c|c|c|c|c|c|c|c|c|c|c|c|c|c|c|}
\hline
k       & 0     &1 &2     &4    &5     &7     &8 & 10    &11    &13    &14    &16 & 17& 19 & 20     & 22  & 23     &25   &26\\
\hline
$|G|\cdot\deg_W$  &32  &0 &22 &12 &34 &24 &46 &36 &58 &48 &16 &6 &28 &18 &40 &30 &52 &42 &64 \\
\hline
invariants & $y^8e_0$& $\1$& $e_2$& $e_4$& $e_5$& $e_7$& $e_8$& $e_{10}$& $e_{11}$& $e_{13}$& $e_{14}$& $e_{16}$& $e_{17}$ & $e_{19}$& $e_{20}$& $e_{22}$& $e_{23}$& $e_{25}$& $e_{26}$\\
\hline
\end{tabular}
\end{center}

Potential non-zero correlators:

By the pairing axiom $\br{\1,\1,e_{26}}$, $\br{\1,e_2,e_{25}}$, $\br{\1,e_4,e_{23}}$, $\br{\1,e_5,e_{22}}$, $\br{\1,e_7,e_{20}}$, $\br{\1,e_8,e_{19}}$, $\br{\1,e_{10},e_{17}}$, $\br{\1,e_{11},e_{16}}$, and $\br{\1,e_{13},e_{14}}$ are equal to 1, and $\br{y^8e_0,y^8e_0,\1}=-1/9$,

By the Concavity axiom $\br{e_2,e_4,e_{22}}$, $\br{e_2,e_7,e_{19}}$, $\br{e_2,e_{10},e_{16}}$, $\br{e_4,e_4,e_{20}}$, $\br{e_4,e_{5},e_{19}}$, $\br{e_4,e_{7},e_{17}}$, $\br{e_4,e_8,e_{16}}$, $\br{e_4,e_{10},e_{14}}$, $\br{e_5,e_7,e_{16}}$, $\br{e_7,e_7,e_{14}}$, $\br{e_{14},e_{16},e_{25}}$, $\br{e_{14},e_{19},e_{22}}$, $\br{e_{16},e_{16},e_{23}}$, $\br{e_{16},e_{17},e_{22}}$, $\br{e_{16},e_{19},e_{20}}$, and $\br{e_{17},e_{19},e_{19}}$ are all equal to 1.

By the composition axiom $\br{y^8e_{0},e_{4},e_{24}}=\pm 1$.

Consider the map $\varphi:\C[X,Y]\rightarrow \h{Z_{19}^T}{G}$ defined by $X\mapsto e_{14}$ and $Y\mapsto e_{16}$ and extending to a $\C$-algebra homomorphism. One can check directly that this map is surjective. 

Our correlators tell us that $e_{14}^3=-9e_{13}$, $e_{16}^8=e_{13}$ and $e_{14}^2e_{16}=0$. 

Hence $(3X^2Y,X^3+9Y^8)\subset\ker\varphi$. Since $\C[X,Y]/(3X^2Y,X^3+9Y^8)=\q{Z_{19}}$ has dimension 19, we deduce that the inclusion is equality, and so we have the isomorphism $\h{Z_{19}^T}{G}\cong \q{Z_{19}}$.
\newline\hrule

\subsubsection*{$\mathbf{W_{17}=x^4+xy^5}$}
\[
\jac=(4x^3+y^5,5xy^4)
\phantom{XXX}
q_x=\tfrac{5}{20}, q_y=\tfrac{3}{20}, \hat c=\tfrac {24}{20}
\phantom{XXX}
G=\br J\cong \Z/20\Z
\]

\[
\Fix{J^k}=
\begin{cases}
\C^2       & \text{if }k=0\\
\C_x        &\text{if }4|k\\
0           &\text{otherwise}
\end{cases}
\phantom{XXX}
\q{}|_{\Fix J^k}=
\begin{cases}
\br{1, x, x^2, y,\dots, y^8, xy, xy^2, xy^3, x^2y, x^2y^2, x^2y^3},\mu=17\\
\br{1,x,x^2},\mu=3\\
\br 1                \\
\end{cases}
\]

\begin{center}
\begin{tabular}{|c||c|c|c|c|c|c|c|c|c|c|c|c|c|c|c|c|}
\hline
$k$        &0&1&2&3&5&6&7&9&10&11&13&14&15&17&18&19\\
\hline
$|G|\cdot\deg_W$      &24 &0 &16 &32 &24 &40 &16 &8&24 &40 &32 &8 &24 &16 &32 &48\\
\hline
invariants & $y^4e_0$& $\1$& $e_{2}$& $e_{3}$& $e_{5}$& $e_{6}$& $e_{7}$& $e_{9}$& $e_{10}$& $e_{11}$& $e_{13}$& $e_{14}$& $e_{15}$& $e_{17}$& $e_{18}$& $e_{19}$\\
\hline
\end{tabular}

\end{center}

Potential non-zero correlators:

By the pairing axiom $\br{\1,\1,e_{19}}$, $\br{\1,e_{2},e_{18}}$, $\br{\1,e_{3},e_{17}}$, $\br{\1,e_{5},e_{15}}$, $\br{\1,e_{6},e_{14}}$,
$\br{\1,e_{7},e_{13}}$, $\br{\1,e_{9},e_{11}}$, and $\br{\1,e_{10},e_{10}}$ are all equal to 1, and $\br{y^4e_{0},y^4e_{0},\1}=-1/5$

By the concavity axiom $\br{e_{2},e_{2},e_{17}}$,
$\br{e_{2},e_{5},e_{14}}$, $\br{e_{2},e_{9},e_{10}}$,
$\br{e_{3},e_{9},e_{9}}$, $\br{e_{5},e_{7},e_{9}}$,
$\br{e_{7},e_{17},e_{17}}$, $\br{e_{9},e_{14},e_{18}}$,
$\br{e_{9},e_{15},e_{17}}$, $\br{e_{10},e_{14},e_{17}}$,
$\br{e_{13},e_{14},e_{14}}$, are all equal to 1.

By the Index Zero axiom $\br{e_{7},e_{7},e_{7}}=-5$.

By the Composition axiom $\br{y^4e_{0},e_{7},e_{14}}=\pm 1$

Consider the map $\varphi:\C[X,Y]\rightarrow \h{W_{17}}{G}$ defined by $X\mapsto e_{14}$ and $Y\mapsto e_{9}$ and extending to a $\C$-algebra homomorphism. One can check directly that this map is surjective. 

Our correlators tell us that $e_{14}^3e_{9}=0$, $e_{14}^4=-5e_{13}$, and $e_{9}^4=e_{13}$. 

Hence $(4X^3Y,X^4+5Y^4)\subset\ker\varphi$. Since $\C[X,Y]/(4X^3Y,X^4+5Y^4)=\q{W_{17}^T}$ has dimension 16, we deduce that the inclusion is equality, and so we have the isomorphism $\h{W_{17}}{G}\cong\q{W_{17}^T}$.
\newline\hrule

\subsubsection*{$\mathbf{W_{17}^{T}=x^4y+y^5}$}
\[
\jac=(4x^3y,x^4+5y^4)
\phantom{XXX}
q_x=\tfrac{4}{20}, q_y=\tfrac{4}{20}, \hat c=\tfrac {24}{20}
\phantom{XXX}
G\cong \Z/20\Z
\phantom{XXX}
\br J\cong \Z/5\Z
\]

Since $J$ does not generate $G$ we will use the generator $\g=(\zeta,\zeta^{-4})$, where $\zeta$ is a primitive 20-th root of unity.
\[
\Fix{\g^k}=
\begin{cases}
\C^2       & \text{if }k=0\\
\C_y        &\text{if }5|k\\
0           &\text{otherwise}
\end{cases}
\phantom{XXX}
\q{}|_{\Fix \g^k}=
\begin{cases}
\br{1, x, x^2, x^3, y, \dots, y^4, xy, \dots, xy^4, x^2y, \dots, x^2y^4},\mu=16 \\
\br{1,y,y^2,y^3},\mu=4 \\
\br 1                     \\
\end{cases}
\]

\begin{center}
\begin{tabular}{|c||c|c|c|c|c|c|c|c|c|c|c|c|c|c|c|c|c|}
\hline
$k$          &0   &1    &2   &3    &4 &6   &7     &8 &9   &11    &12  &13    &14  &16   &17    &18  &19\\
\hline
$|G|\cdot\deg_W$      &24 &18 &12 &6 &0 &28 &22 &16 &10 &38 &32 &26 &20   &48 &42 &36 &30\\
\hline
invariants & $x^3e_0$& $e_{1}$& $e_{2}$& $e_{3}$& $\1$& $e_{6}$& $e_{7}$& $e_{8}$& $e_{9}$& $e_{11}$& $e_{12}$& $e_{13}$& $e_{14}$& $e_{16}$& $e_{17}$& $e_{18}$& $e_{19}$\\
\hline
\end{tabular}
\end{center}

Potential non-zero correlators:

By the pairing axiom $\br{\1,e_{1},e_{19}}$, $\br{\1,e_{2},e_{18}}$, $\br{\1,e_{3},e_{17}}$, $\br{\1,\1,e_{16}}$, $\br{\1,e_{6},e_{14}}$, $\br{\1,e_{7},e_{13}}$, $\br{\1,e_{8},e_{12}}$, and $\br{\1,e_{9},e_{11}}$ are all equal to 1, and $\br{x^3e_{0},x^3e_{0},\1}=-1/4$

By the concavity axiom $\br{e_{1},e_{9},e_{14}}$, $\br{e_{2},e_{3},e_{19}}$, $\br{e_{2},e_{8},e_{14}}$, $\br{e_{2},e_{9},e_{13}}$, $\br{e_{3},e_{3},e_{18}}$, $\br{e_{3},e_{7},e_{14}}$, $\br{e_{3},e_{8},e_{13}}$, $\br{e_{3},e_{9},e_{12}}$, $\br{e_{6},e_{9},e_{9}}$, $\br{e_{7},e_{8},e_{9}}$, and $\br{e_{8},e_{8},e_{8}}$ are all equal to 1.

By the Index Zero axiom $\br{e_{1},e_{1},e_{2}}=-4$.

By the composition axiom
$\br{x^3e_{0},e_{1}=e_{3}}\br{x^3e_{0},e_{2},e_{2}}=\pm 1$

Consider the map $\varphi:\C[X,Y]\rightarrow \h{W_{17}^T}{G}$ defined by $X\mapsto e_{9}$ and $Y\mapsto e_{3}$ and extending to a $\C$-algebra homomorphism. One can check directly that this map is surjective. 

Our correlators tell us that $e_{9}e_{3}^4=0$, $e_{3}^5=-4e_{19}$, and $e_{9}^3=e_{19}$. 

Hence $(4X^3+Y^5,5XY^4)\subset\ker\varphi$. Since $\C[X,Y]/(4X^3+Y^5,5XY^4)=\q{}$ has dimension 17, we deduce that the inclusion is equality, and so we have the isomorphism $\h{W_{17}^T}{G}\cong\q{W_{17}}$.
\newline\hrule

\subsubsection*{$\mathbf{Q_{17}=x^3+xy^5+yz^2}$}
\[
\jac=(3x^2+y^5,5xy^4+z^2,2yz)
\phantom{XXX}
q_x=\tfrac{10}{30}, q_y=\tfrac{4}{30}, q_z=\tfrac{13}{30}, \hat c=\tfrac {36}{30}
\phantom{XXX}
G=\br J\cong \Z/30\Z
\]
\[
\Fix{J^k}=
\begin{cases}
\C^3       & \text{if }k=0\\
\C^2_{xy}   & \text{if }k=15\\
\C_x        &\text{if }3|k,k\neq 15\\
0           &\text{otherwise}
\end{cases}
\phantom{XXX}
\q{}|_{\Fix J^k}=
\begin{cases}
\br{1, x, z, z^2, y, y^2, \dots, y^9, xy, xy^2, xy^3, xz},\mu=17\\
\br{1, x, y, y^2, \dots, y^8, xy, xy^2, xy^3},\mu=13 \\
\br{1,x},\mu=2 \\
\br 1             \\
\end{cases}
\]

\begin{center}
\begin{tabular}{|c||c|c|c|c|c|c|c|c|c|c|c|}
\hline
$k$        &1 &2 &4 &5 &7 &8 &10 &11 &13 &14 &15\\
\hline
$|G|\cdot\deg_W$     & 0& 54& 42& 36& 24& 18& 6 &60 &48 &42 &36\\
\hline
invariants & $\1$& $e_{2}$& $e_{4}$& $e_{5}$& $e_{7}$& $e_{8}$& $e_{10}$& $e_{11}$& $e_{13}$& $e_{14}$& $y^4e_{15}$\\
\hline
\end{tabular}

\vspace*{\baselineskip}

\begin{tabular}{|c||c|c|c|c|c|c|c|c|c|c|c}
\hline
$k$        &16&17&19&20&22&23&25&26&28&29 \\
\hline
$|G|\cdot\deg_W$      &30 &24 &12 &66 &54 &48 &36 &30 &18 &72\\
\hline
invariants & $e_{16}$& $e_{17}$& $e_{19}$& $e_{20}$& $e_{22}$& $e_{23}$& $e_{25}$& $e_{26}$& $e_{28}$& $e_{29}$\\
\hline
\end{tabular}
\end{center}

Potential non-zero correlators:

By the pairing axiom $\br{\1,\1,e_{29}}$, $\br{\1,e_{2},e_{28}}$, $\br{\1,e_{4},e_{26}}$, $\br{\1,e_{5},e_{25}}$, $\br{\1,e_{7},e_{23}}$,
$\br{\1,e_{8},e_{22}}$, $\br{\1,e_{10},e_{20}}$, $\br{\1,e_{11},e_{19}}$, $\br{\1,e_{13},e_{17}}$, and $\br{\1,e_{14},e_{16}}$ are all equal to 1, and $\br{y^4e_{15},y^4e_{15},\1}=-1/5$

By the concavity axiom $\br{e_{2},e_{10},e_{19}}$,
$\br{e_{4},e_{8},e_{19}}$, $\br{e_{4},e_{10},e_{17}}$,
$\br{e_{5},e_{10},e_{16}}$, $\br{e_{7},e_{8},e_{16}}$,
$\br{e_{8},e_{10},e_{13}}$, $\br{e_{8},e_{25},e_{28}}$,
$\br{e_{10},e_{10},e_{11}}$, $\br{e_{10},e_{23},e_{28}}$,
$\br{e_{10},e_{25},e_{26}}$, $\br{e_{16},e_{17},e_{28}}$,
$\br{e_{16},e_{19},e_{26}}$, $\br{e_{17},e_{19},e_{25}}$,
$\br{e_{19},e_{19},e_{23}}$ are all equal to 1.

By the index zero axiom $\br{e_{5},e_{7},e_{19}}$,
$\br{e_{5},e_{28},e_{28}}$, $\br{e_{7},e_{7},e_{17}}$,
$\br{e_{7},e_{10},e_{14}}$, $\br{e_{7},e_{26},e_{28}}$,
$\br{e_{14},e_{19},e_{28}}$, are all equal to $-2$.

By the composition axiom $\br{e_{8},e_{8},y^4e_{15}}=\pm 1$.

Consider the map $\varphi:\C[X,Y,Z]\rightarrow \h{Q_{17}}{G}$ defined by $X\mapsto e_{8}$ and $Y\mapsto e_{10}$ and $Z\mapsto e_{16}$ and extending to a $\C$-algebra homomorphism. One can check directly that this map is surjective. 

Our correlators tell us that $e_{8}^2e_{10}=0$, $e_{10}^5=-2e_{16}$, $e_8^3=-5e_{22}$, and $e_{10}^4e_{16}=e_{22}$.

Hence $(3X^2,X^3+5Y^4Z,Y^5+2Z)\subset\ker\varphi$. Since $\C[X,Y]/(3X^2,X^3+5Y^4Z,Y^5+2Z)=\q{Q_{17}^T}$ has dimension 21, we deduce that the inclusion is equality, and so we have the isomorphism $\h{Q_{17}}{G}\cong\q{Q_{17}^T}$.
\newline\hrule

\subsubsection*{$\mathbf{Q_{17}^{T}=x^3y+y^5z+z^2}$}
\[
\jac=(3x^2y,x^3+5y^4z,y^5+2z)
\phantom{XXX}
q_x=\tfrac{9}{30}, q_y=\tfrac{3}{30}, q_z=\tfrac{15}{30}, \hat c=\tfrac {36}{30}
\phantom{XXX}
G\cong \Z/30\Z
\phantom{XXX}
\br J\cong \Z/10\Z
\]

Since $J$ does not generate $G$ we will use the generator $\g=(\zeta, \zeta^{-3}, \zeta^{15})$, where $\zeta$ is a primitive $30$-th root of unity.

\[
\Fix{\g^k}=
\begin{cases}
\C^3       & \text{if }k=0\\
\C_z        &\text{if }2|k,10\nmid k\\
\C^2_{yz}         & \text{if }10|k\\
0           &\text{otherwise}
\end{cases}
\phantom{XXX}
\q{}|_{\Fix \g^k}=
\begin{cases}
\br{1, x, x^2, z, y, \dots, y^4, xy, \dots, xy^4, xz, yz, \dots, y^4z, xyz, \dots, xy^4z}\\
\br{1} \\
\br{1, y, z, y^2, y^3, y^4, yz, y^2z, y^3z} ,\mu=9\\
\br 1                    \\
\end{cases}
\]

\begin{center}
\begin{tabular}{|c||c|c|c|c|c|c|c|c|c|c|c|c|c|c|c|c|c|}
\hline
$k$        &1&3&5&7&9&10&11&13&15&17&19&20&21&23&25&27&29\\
\hline
$|G|\cdot\deg_W$     &32 &24 &16 &8 &0 &26 &52 &44&36 &28 &20 &46 &72 &64 &56 &48 &40\\
\hline
invariants& $e_1$& $e_{3}$& $e_{5}$& $e_{7}$& $\1$& $y^4e_{10}$& $e_{11}$& $e_{13}$& $e_{15}$& $e_{17}$& $e_{19}$& $y^4e_{20}$& $e_{21}$& $e_{23}$& $e_{25}$& $e_{27}$& $e_{29}$\\
\hline
\end{tabular}

\end{center}

Potential non-zero correlators:

By the pairing axiom $\br{\1,e_{1},e_{29}}$, $\br{\1,e_{3},e_{27}}$, $\br{\1,e_{5},e_{25}}$, $\br{\1,e_{7},e_{23}}$, $\br{\1,\1,e_{21}}$,
$\br{\1,e_{11},e_{19}}$, $\br{\1,e_{13},e_{17}}$, and $\br{\1,e_{15},e_{15}}$ are equal to 1, and $\br{y^4e_{10},y^4e_{20},\1}=-1/5$.

By the concavity axiom we see that  $\br{e_{1},e_{19},e_{19}}$, $\br{e_{3},e_{7},e_{29}}$, $\br{e_{3},e_{17},e_{19}}$, $\br{e_{5},e_{5},e_{29}}$, $\br{e_{5},e_{7},e_{27}}$, $\br{e_{5},e_{15},e_{19}}$, $\br{e_{5},e_{17},e_{17}}$, $\br{e_{7},e_{7},e_{25}}$, $\br{e_{7},e_{13},e_{19}}$, and $\br{e_{7},e_{15},e_{17}}$ are all equal to 1.

By the Index Zero axiom $\br{e_{1},e_{1},e_{7}}=\br{e_{1},e_{3},e_{5}}=\br{e_{3},e_{3},e_{3}}=-3$.

We cannot compute the last three-point correlator $a=\br{y^4e_{10},y^4e_{10},e_{19}}$. If $a\neq 0$, then we can establish an isomorphism between the graded rings $\h{Q_{17}^T}{G}$ and $\q{Q_{17}}$ in the following way:

Consider the map $\varphi:\C[X,Y,Z]\rightarrow \h{Q_{17}^T}{G}$ defined by $X\mapsto e_{19}$, $Y\mapsto e_{7}$ and $Z\mapsto \beta y^4e_{10}$ with $\beta^2=\frac {-5}{a}$. Extending this map by $\C$-linearity and multiplicativity we get a $\C$-algebra homomorphism.One can check directly that this map is surjective. 

Our correlators tell us that $e_{19}^2=e_{29}$, $e_{7}^5=-3e_{29}$, $e_{19}e_{7}^4=e_{11}$, $(\beta y^4e_{10})^2=\beta^2(y^4e_{10})^2=(\beta^2a) e_{11}=-5e_{11}$, and $e_{7}(y^4e_{10})=0$. 

Hence $(3X^2+Y^5,5XY^4+Z^2,2YZ)\subset\ker\varphi$. Since $\C[X,Y,Z]/(3X^2+Y^5,5XY^4+Z^2,2YZ)=\q{Q_{17}}$ has dimension 17, we deduce that the inclusion is equality, and so we have the isomorphism $\h{Q_{17}^T}{G}\cong\q{Q_{17}}$ under the assumption that $a\neq 0$. 
\newline\hrule

\subsubsection*{$\mathbf{S_{16}=x^2z+yz^2+xy^4}$}
\[
\jac=(2xz+y^4,z^2+4xy^3,x^2+2yz)
\phantom{XXX}
q_x=\tfrac{5}{17}, q_y=\tfrac{3}{17}, q_z=\tfrac{7}{17}, \hat c=\tfrac{21}{17}
\phantom{XXX}
G=\br J\cong \Z/17\Z
\]
\[
\Fix{J^k}=
\begin{cases}
\C^3       & \text{if }k=0\\
0           &\text{otherwise}
\end{cases}
\phantom{XXX}
\q{}|_{\Fix J^0}=
\begin{cases}
\br{1, x, y, y^2, \dots, y^6, z, z^2, z^3, xy, xy^2, yz, y^2z, y^3z}, \mu=16 \\
\br 1           \\
\end{cases}
\]

\begin{center}
\begin{tabular}{|c||c|c|c|c|c|c|c|c|c|c|c|c|c|c|c|c|}
\hline
$k$        &1&2&3&4&5&6&7&8&9&10&11&12&13&14&15&16\\
\hline
$|G|\cdot\deg_W$     &0&30 &26 &22 &18 &14 &10 &6 &36 &32 &28 &24 &20 &16 &12 &42\\
\hline
invariants& $\1$& $e_{2}$& $e_{3}$& $e_4$& $e_{5}$& $e_{6}$& $e_{7}$& $e_8$& $e_{9}$& $e_{10}$& $e_{11}$& $e_{12}$& $e_{13}$& $e_{14}$& $e_{15}$& $e_{16}$\\
\hline
\end{tabular}

\end{center}

Potential non-zero correlators:

By the pairing axiom $\br{\1,\1,e_{16}}$, $\br{\1,e_{2},e_{15}}$, $\br{\1,e_{3},e_{14}}$, $\br{\1,e_{4},e_{13}}$, $\br{\1,e_{5},e_{12}}$,
$\br{\1,e_{6},e_{11}}$, $\br{\1,e_{7},e_{10}}$, and $\br{\1,e_{8},e_{9}}$ are equal to 1.

By the concavity axiom $\br{e_{2},e_{8},e_{8}}$, $\br{e_{3},e_{7},e_{8}}$, $\br{e_{4},e_{6},e_{8}}$, $\br{e_{5},e_{6},e_{7}}$, $\br{e_{6},e_{14},e_{15}}$, $\br{e_{7},e_{13},e_{15}}$, $\br{e_{8},e_{12},e_{15}}$, and $\br{e_{8},e_{13},e_{14}}$ are all equal to 1.

By the Index Zero axiom
$\br{e_{4},e_{7},e_{7}}$, $\br{e_{7},e_{14},e_{14}}$,
$\br{e_{5},e_{5},e_{8}}$, $\br{e_{5},e_{15},e_{15}}$ all equal to
-2, and $\br{e_{6},e_{6},e_{6}}=-4$.

Consider the map $\varphi:\C[X,Y,Z]\rightarrow \h{S_{16}}{G}$ defined by $X\mapsto e_{7}$, $Y\mapsto e_{8}$ and $Z\mapsto e_{6}$ and extending to a $\C$-algebra homomorphism. One can check directly that this map is surjective. 

Our correlators tell us that $e_{7}e_{6}=e_{12}$, $e_{8}^4=-2e_{12}$, $e_{6}^2=-4e_{11}$, $e_{7}e_{8}^3=e_{11}$,
$e_{7}^2=-2e_{13}$, and $e_{8}e_{6}=e_{13}$. 

Hence $(2XZ+Y^4,Z^2+4ZY^3,X^2+2YZ)\subset\ker\varphi$. Since $\C[X,Y,Z]/(2XZ+Y^4,Z^2+4ZY^3,X^2+2YZ)=\q{S_{16}}$ has dimension 16, we deduce that the inclusion is equality, and so we have the isomorphism $\h{S_{16}}{G}\cong\q{S_{16}}$.
\newline\hrule

\subsubsection*{$\mathbf{S_{17}=x^2z+yz^2+y^6}$}
\[
\jac=(2xy, z^2+6y^5 , x^2+2yz)
\phantom{XXX}
q_x=\tfrac{7}{24}, q_y=\tfrac{4}{24}, q_z=\tfrac{10}{24}, \hat{c}= \tfrac{30}{24}
\phantom{XXX}
G=\br J\cong \Z/24\Z
\]
\[
\Fix{J^k}=\begin{cases}
\C^3  & \text{if}\quad k=0\\
\C^2_{yz} & \text{if}\quad k= 12\\
\C_y  & \text{if}\quad 6|k,k\neq 12\\
0  & \text{otherwise}
\end{cases}
\phantom{XXX}
\q{}|_{\Fix J^k}=
\begin{cases}
\br{1,x,y,\dots,y^4,z,z^2,z^3,xy,\dots,xy^4,yz,\dots, y^4 z },\mu=17 \\
\br{1,y,y^2 , y^3 , y^4 ,z, z^2},\mu=7 \\
\br{1,y, y^2 , y^3 , y^4},\mu=5 \\
\br 1 \\
\end{cases}
\]

\begin{center}
\begin{tabular}{|c||c|c|c|c|c|c|c|c|c|c|c|c|c|c|c|c|c|c|c|c|c|c|}
\hline
$k$        &  1 & 2 & 3 & 4 & 5 &  7 & 8 & 9 & 10 & 11\\
\hline
$|G|\cdot\deg_W$    &0 &42 &36 &30 &24   &12 &6 &48 &42 &36\\
\hline
invariants& $\1$& $e_2$& $e_3$& $e_4$& $e_5$& $e_7$& $e_8$& $e_9$& $e_{10}$& $e_{11}$\\
\hline
\end{tabular}

\vspace*{\baselineskip}

\begin{tabular}{|c||c|c|c|c|c|c|c|c|c|c|c|c|c|c|c|c|c|c|c|c|c|c|c|}
\hline
$k$        & 12 & 13 & 14 & 15 & 16 & 17 & 19 & 20 & 21 & 22 & 23  \\
\hline
$|G|\cdot\deg_W$    &30  &24   &18 &12 &54 &48 &36 &30 &24   &42 &60\\
\hline
invariants & $ze_{12}$& $e_{13}$& $e_{14}$& $e_{15}$& $e_{16}$& $e_{17}$& $e_{19}$& $e_{20}$& $e_{21}$& $e_{22}$& $e_{23}$\\
\hline
\end{tabular}
\end{center}

Potential non-zero correlators:

By the pairing axiom $\br{\1,\1,e_{23}}$, $\br{\1,e_2,e_{22}}$, $\br{\1,e_3,e_{21}}$, $\br{\1,e_4,e_{20}}$, $\br{\1,e_5,e_{19}}$,  $\br{\1,e_7,e_{17}}$, $\br{\1,e_8,e_{16}}$, $\br{\1,e_{9},e_{15}}$, $\br{\1,e_{10},e_{14}}$, and $\br{\1,e_{11},e_{13}}$ are equal to 1, and $\br{ze_12,ze_12,\1}=-1/2$,

By the Concavity axiom $\br{e_2,e_8,e_{15}}$, $\br{e_3,e_7,e_{15}}$, $\br{e_3,e_{8},e_{14}}$, $\br{e_4,e_8,e_{13}}$, $\br{e_5,e_{7},e_{13}}$, $\br{e_7,e_{8},e_{10}}$, $\br{e_7,e_{20},e_{22}}$, $\br{e_8,e_{8},e_{9}}$, $\br{e_8,e_{19},e_{22}}$, $\br{e_8,e_{20},e_{21}}$, $\br{e_{13},e_{14},e_{22}}$, $\br{e_{13},e_{15},e_{21}}$, $\br{e_{14},e_{15},e_{20}}$, $\br{e_{15},e_{15},e_{19}}$ are all equal to 1.

By the index zero axiom, $\br{e_4,e_7,e_{14}}$, $\br{e_5,e_5,e_{15}}$, $\br{e_5,e_{22}, e_{22}}$, $\br{e_7,e_7,e_{11}}$, $\br{e_7,e_{21},e_{21}}$, and$\br{e_{14},e_{14},e_{21}}$  to be equal to -2

By the composition axiom $\br{e_{5},e_{8},ze_{12}}=\br{ze_{12},e_{15},e_{22}}=\pm 1$.

Consider the map $\varphi:\C[X,Y,Z]\rightarrow \h{S_{17}}{G}$ defined by $X\mapsto e_8$, $Y\mapsto e_{13}$, and $Z\mapsto e_7$ and extending to a $\C$-algebra homomorphism. One can check directly that this map is surjective. 

Our correlators tell us that $e_{13}*e_7=e_{19}$ and $e_8^6=-2e_{19}$, $e_7^2=-2e_{13} $ and $e_8^5*e_7=0$. 

Hence $(6X^5Z, Z^2+2Y,X^6+2YZ) \subset \ker \varphi$. Since $\C[X,Y,Z]/(6X^5Z, Z^2+2Y, X^6+2YZ)=\q{S_{17}^T}$ has dimension 21, we deduce that the inclusion is equality, and so we have the isomorphism $\h{S_{17}}{G}\cong\q{S_{17}^T}$.
\newline\hrule

\subsubsection*{$\mathbf{S_{17}^T=x^6z+z^2y+y^2}$}
\[
\jac=(6x^5z,z^2+2y,x^6+2zy)
\phantom{XXX}
q_x=\tfrac{3}{24}, q_y=\tfrac{12}{24}, q_z=\tfrac{6}{24}, \hat c=\tfrac {30}{24}
\phantom{XXX}
G\cong \Z/24\Z
\phantom{XXX}
\br J\cong \Z/8\Z
\]

Since $J$ does not generate the symmetry group, we will use the generator $\g=(\zeta, \zeta^{12},\zeta^{-6})$ where $\zeta$ is a primitive $24$-th root of unity.
\[
\Fix{\g^k}=\begin{cases}
\C^3 & \text{if}\quad k=0\\
\C_y & \text{if}\quad 2\mid k, 4\nmid k\\
\C^2_{yz} & \text{if}\quad 4\mid k\\
0    & \text{otherwise}
\end{cases}
\phantom{XXX}
\q{}|_{\Fix \g^k}=
\begin{cases}
\br{1,x,x^2,\dots,x^{10},z,z^2,xz,x^2z,x^3z,x^4z,xz^2,\dots,x^4z^2} \\
\br{1} \\
\br{1,z,z^2} \\
\br 1  \\
\end{cases}
\]

\begin{center}
\begin{tabular}{|c||c|c|c|c|c|c|c|c|c|c|c|c|c|c|c|c|c|c|c|c|c|c|c|}
\hline
$k$        & 1    & 3 & 4    & 5    & 7   & 8     & 9   & 11  & 12 & 13    & 15  & 16    & 17    & 19  & 20    & 21  & 23\\
\hline
$|G|\cdot\deg_W$    & 20  & 0 & 14 & 28  & 8 & 22 & 36 & 16 & 30 & 44  & 24   & 38 & 52  & 32 & 46 & 60 & 40 \\
\hline
invariants & $e_1$& $\1$& $ze_4$& $e_5$& $e_7$& $ze_8$& $e_9$& $e_{11}$& $ze_{12}$& $e_{13}$& $e_{15}$& $ze_{16}$& $e_{17}$& $e_{19}$& $ze_{20}$& $e_{21}$& $e_{23}$\\
\hline
\end{tabular}
\end{center}

Potential non-zero correlators:

By the Pairing axiom $\br{\1,e_1,e_{23}}$, $\br{\1,\1,e_{21}}$, $\br{\1,e_{5},e_{19}}$, $\br{\1,e_7,e_{17}}$, $\br{\1,e_{9},e_{15}}$, and $\br{\1,e_{11},e_{13}}$ are equal to one.
$\br{ze_4,ze_{20},\1}$, $\br{ze_8,ze_{16},\1}$, and
$\br{ze_{12},ze_{12},\1}$ are equal to $-1/2$.

By the Concavity axiom $\br{e_{1},e_{5},e_7}$, $\br{e_{1},e_7,e_{19}}$, $\br{e_{5},e_{7},e_{15}}$, $\br{e_7,e_7,e_{13}}$, $\br{e_7,e_{9},e_{11}}$,  $\br{e_{11},e_{11},e_{5}}$, $\br{e_{11},e_{15},e_{1}}$, are all equal to 1.

By the index zero axiom, $\br{e_1,e_1,e_{1}}$ is equal to -6.
The rest of the three-point correlators we cannot compute. However, if we put $c=\br{ze_4,ze_4,e_{19}}$, $a=\br{ze_4,e_7,ze_{16}}$, $e=\br{ze_4,ze_8,e_{15}}$, $d=\br{ze_4,e_{11},ze_{12}}$, $b=\br{e_7,ze_8,ze_{12}}$, $f=\br{ze_8,ze_8,e_{11}}$ the composition axiom gives us the following relations:
\begin{align}
-2af&=-2bd=e\\
-2ab&=d\\
-2ae&=c\\
-2b^2&=f\\
-2d^2&=c
\end{align}

If $c\neq 0$, then we can construct an isomorphism. Consider the map $\varphi:\C[X,Y,Z]\rightarrow \h{S_{17}^T}{G}$ defined by $X\mapsto \beta ze_4$, where $\beta^2=\frac{-2}{C}$, $Y\mapsto e_7$, and $Z\mapsto e_1$ and extending to a $\C$-algebra homomorphism. One can check directly that this map is surjective. 

Our correlators tell us that $(ze_4)e_1=0$, $e_1^2=-6e_{23}$, $e_7^5=e_{23}$, $(\beta ze_4)^2=\beta^2Ce_{5}=-2e_5$ and $e_7e_1=e_5$. 

Hence $(2XZ, Z^2+6Y^5, X^2+2YZ)\subset\ker\varphi$. Since $\C[X,Y,Z]/(2XZ, Z^2+6Y^5, X^2+2YZ)=\q{S_{17}}$ has dimension 17, we deduce that the inclusion is equality, and so we have the isomorphism $\h{S_{17}^T}{G}\cong\q{S_{17}}$ if $c\neq 0$.

\subsection{Singularities of Corank 3}\label{co3}
In this section, we consider four singularities $Q_{2,0}$ and $S_{1,0}$, and the transpose of each. The other two in Arnol'd's list are not quasi-homogeneous for any choice of constants.

\subsubsection*{$\mathbf{Q_{2,0}=x^3 + yz^2 + xy^4}$}
\[
\jac=(3x^2+y^4,z^2+4xy^3,2yz)
\phantom{XXX}
q_x=\tfrac 8{24}, q_y=\tfrac 4{24}, q_z=\tfrac {10}{24}, \hat c=\tfrac{28}{24}
\phantom{XXX}
G\cong \Z/24\Z
\phantom{XXX}
\br J\cong \Z/12\Z
\]

Since $J$ does not generate the full group of diagonal symmetries we will use
$\g=(\zeta^8,\zeta^{-2},\zeta)$, where $\zeta^{24}=1$.

\[
\Fix {\g^k} = \begin{cases}
 \C^3 & \text{if}\quad k = 0\\
\C^2_{xy} & \text{if}\quad k=12\\
\C_x & \text{if}\quad 3|k, k\neq 12\\
  0    & \text{otherwise}
\end{cases}
\phantom{XXX}
\q{}|_{\Fix \g^k} =
\begin{cases}
\br{1,x,y,y^2\dots, y^7, z,z^2,xy,xz},\mu=14 \\
\br{1,x,x^2,y,y^2,y^3,xy,xy^2,x^2y,x^2y^2},\mu=10 \\
\br{1,x},\mu=2 \\
\br{1} \\
\end{cases}
\]

\begin{center}
\begin{tabular}{|c||c|c|c|c|c|c|c|c|c|c|c|c|c|c|c|c|c|c|c|c|c|c|c|c|}
\hline
k       & 1 & 2 & 4 & 5 & 7 & 8 & 10 & 11 & 12 & 13 & 14 & 16 & 17 & 19 & 20 & 22 & 23\\
\hline
$|G|\cdot\deg_W$  & 18   & 32 & 12 & 26 & 6 & 20 & 0 & 14 & 28 & 42 & 56 & 36  & 50 & 30 & 44 & 24 & 38\\
\hline
invariants & $e_1$& $e_2$& $e_4$& $e_5$& $e_7$& $e_8$& $\1$& $e_{11}$& $y^3e_{12}$& $e_{13}$& $e_{14}$& $e_{16}$& $e_{17}$& $e_{19}$& $e_{20}$& $e_{22}$& $e_{23}$\\
\hline
\end{tabular}
\end{center}

Potential non-zero correlators:

By the pairing axiom $\br{\1,e_{1},e_{23}}$, $\br{\1,e_{2},e_{22}}$, $\br{\1,e_{4},e_{20}}$, $\br{\1,e_5,e_{19}}$, $\br{\1,e_7,e_{17}}$, $\br{\1,e_8,e_{16}}$, $\br{\1,e_{11},e_{13}}$, and $\br{\1,\1,e_{14}}$ are all equal to 1, and $\br{\1,y^3e_{12},y^3e_{12}}=-1/4$.

By the Concavity Axiom $\br{e_{1},e_{11},e_{22}}$, $\br{e_4,e_7,e_{23}}$, $\br{e_4,e_8,e_{22}}$, $\br{e_4,e_{11},e_{19}}$, $\br{e_5,e_{7},e_{22}}$, $\br{e_7,e_7,e_{20}}$, $\br{e_7,e_8,e_{19}}$, $\br{e_7,e_{11},e_{16}}$
and are all equal to 1.

By the Index Zero Axiom $\br{e_1,e_2,e_7}=-2$

By the Composition Axiom $\br{e_{11},e_{11},y^3e_{12}}=\pm 1$.

Consider the map $\varphi:\C[X,Y,Z]\rightarrow \h{Q_{2,0}}{G}$ defined by $X\mapsto e_{11}$, $Y\mapsto e_7$ and $Z\mapsto e_{22}$ and extending to a $\C$-algebra homomorphism. One can check directly that this map is surjective. 

Our correlators tell us that $e_{11}^2e_7= 0$,
$e_{11}^3=-4e_{13}$, $e_7^3e_{22}=e_{13}$,and that $e_7^4=-2e_{22}$. 

Hence $(3X^2Y,X^3+4Y^3Z,2Z+Y^4)\subset\ker\varphi$. Since $\C[X,Y,Z]/(3X^2Y,X^3+4Y^3Z,2Z+Y^4)=\q{Q_{2,0}^T}$ has dimension 17, we deduce that the inclusion is equality, and so we have the isomorphism $\h{Q_{2,0}}{G}\cong \q{Q_{2,0}^T}$.
\newline\hrule

\subsubsection*{$\mathbf{Q_{2,0}^T=x^3y+y^4z+z^2}$}
\[
\jac=(3x^2y,x^3+4y^3z,2z+y^4)
\phantom{XXX}
q_x=\tfrac 7{24}, q_y=\tfrac 3{24}, q_z=\tfrac {12}{24}, \hat c=\tfrac {28}{24}
\phantom{XXX}
G =\br J\cong\Z/24\Z
\]
\[
\Fix{J^k}=\begin{cases}
\C^3 & \text{if}\quad k=0\\
\C^2_{yz} & \text{if}\quad 8|k\\
\C_z  & \text{if}\quad 2|k, 8\nmid k\\
0    & \text{otherwise}
\end{cases}
\phantom{XXX}
\q{}|_{\Fix J^k}=
\begin{cases}
\br{1,x,x^2,y,y^2,\dots,y^7,xy,xy^2,\dots,xy^7},\mu=17 \\
\br{1} \\
\br{1,y,y^2,y^3,z,yz,y^2z},\mu=7 \\
\br 1  \\
\end{cases}
\]

\begin{center}
\begin{tabular}{|c||c|c|c|c|c|c|c|c|c|c|c|c|c|c|c|c|c|c|c|c|c|c|c|c|}
\hline
$k$        & 1 & 3 & 5 & 7 & 8 & 9 & 11 & 13 & 15 & 16 & 17 & 19 & 21 & 23\\
\hline
$|G|\cdot\deg_W$    & 0 & 40 & 32& 24 & 20 & 16 & 8& 48 & 40 & 36 & 32 & 24 & 16 & 56\\
\hline
invariants& $\1$& $e_3$& $e_5$& $e_7$& $y^3e_8$& $e_9$& $e_{11}$& $e_{13}$& $e_{15}$& $y^3e_{16}$& $e_{17}$& $e_{19}$& $e_{21}$& $e_{23}$\\
\hline
\end{tabular}
\end{center}

Potential non-zero correlators:

By the paring axiom $\br{\1,\1,e_{23}}$, $\br{\1,e_3,e_{21}}$,
$\br{\1,e_5,e_{19}}$, $\br{\1,e_7,e_{17}}$, $\br{\1,e_9,e_{15}}$, and
$\br{\1,e_{11},e_{13}}$ are equal to 1, and $\br{y^3e_8,y^3e_{16},\1}=-1/4$.

By the concavity axiom $\br{e_3,e_{11},e_{11}}$, $\br{e_5,e_9,e_{11}}$, $\br{e_7,e_9,e_9}$, $\br{e_9,e_{19},e_{21}}$, $\br{e_{11},e_{17},e_{21}}$, and $\br{e_{11},e_{19},e_{19}}$ are equal to 1.

By the Index Zero axiom $\br{e_7,e_7,e_{11}}=\br{e_{7},e_{21},e_{21}}=-3$

The axioms do not provide us an easy way to calculate the last three-point correlator, $a=\br{y^3e_8,y^3e_8,e_9}$. If $a\neq 0$, then we can construct the desired homomorphism as follows:

Consider the map $\varphi:\C[X,Y,Z]\rightarrow \h{Q_{2,0}^T}{G}$ defined by $X\mapsto e_9$, where $\beta^2=\frac{-4}a$, $Y\mapsto e_{11}$ and $Z\mapsto \beta y^3e_8$ and extending to a $\C$-algebra homomorphism. One can check directly that this map is surjective. 

Our correlators tell us that $e_{11}(y^3e_8)= 0$,
$e_9^2=e_{17}$, $e_{11}^4= -3e_{17}$, $(\beta y^3e_8)^2=(\beta^2a)e_{15}=-4e_{15}$ and that $e_9e_{11}^3=e_{15}$. 

Hence $(3X^2+Y^4,2YZ,Z^2+4XY^3)\subset\ker\varphi$. Since $\C[X,Y,Z]/(3X^2+Y^4,2YZ,Z^2+4XY^3)=\q{Q_{2,0}}$ has dimension 14, we deduce that the inclusion is equality, and so we have the isomorphism $\h{Q_{2,0}^T}{G}\cong\q{Q_{2,0}}$, assuming $a\neq 0$.
\newline\hrule

\subsubsection*{$\mathbf{S_{1,0}=x^2z + yz^2 + y^5}$}
\[
\jac=(2xz,z^2+5y^4,x^2+2yz)
\phantom{XXX}
q_x=\tfrac 6{20}, q_y=\tfrac 4{20}, q_z=\tfrac 8{20}, \hat c=\tfrac{24}{20}
\phantom{XXX}
G\cong \Z/20\Z
\phantom{XXX}
\br J\cong \Z/10\Z
\]

Since $J$ does not generate the full group of diagonal symmetries we will use $\g=(\zeta,\zeta^4,\zeta^{-2})$, where $\zeta^{20}=1$.
\[
\Fix {\g^k} = \begin{cases}
\C^3 & \text{if}\quad k = 0\\
\C^2_{yz} & \text{if}\quad k=10\\
\C_y & \text{if}\quad 5|k, k\neq 10\\
   0    & \text{otherwise}
\end{cases}
\phantom{XXX}
\q{}|_{\Fix \g^k} =
\begin{cases}
\br{1,x,y,y^2,y^3,z,z^2,z^3,xy,xy^2,y^3,yz,y^2z,y^3z},\mu=14 \\
\br{1,y,y^2,y^3,y^4,z},\mu=6 \\
\br{1,y,y^2,y^3},\mu=4 \\
\br{1} \\
\end{cases}
\]

\begin{center}
\begin{tabular}{|c||c|c|c|c|c|c|c|c|c|c|c|c|c|c|c|c|c|c|c|c|c|c|c|c|}
\hline
k       & 1 & 2 & 3 & 4 & 6 & 7 & 8 & 9 & 10& 11 & 12 & 13 & 14 & 16 & 17 & 18 & 19\\
\hline
$|G|\cdot\deg_W$  & 10 & 16 & 22 & 28 & 0 & 6 & 12 & 18 & 24& 30 & 36 & 42 & 48 & 20 & 26 & 32 & 38\\
\hline
invariants& $e_1$& $e_2$& $e_3$& $e_4$& $\1$& $e_7$& $e_8$& $e_9$& $ze_{10}$& $e_{11}$& $e_{12}$& $e_{13}$& $e_{14}$& $e_{16}$& $e_{17}$& $e_{18}$& $e_{19}$\\
\hline
\end{tabular}
\end{center}

Potential non-zero correlators:

By the Pairing Axiom $\br{\1,e_1,e_{19}}$, $\br{\1,e_2,e_{18}}$, $\br{\1,e_3,e_{17}}$, $\br{\1,e_4,e_{16}}$, $\br{\1,\1,e_{14}}$, $\br{\1,e_7,e_{13}}$, $\br{\1,e_8,e_{12}}$, and $\br{\1,e_9,e_{11}}$ are equal to 1, and $\br{\1,ze_{10},ze_{10}}=-1/2$.

By the Concavity Axiom $\br{e_1,e_7,e_{18}}$, $\br{e_1,e_8,e_{17}}$, $\br{e_1,e_9,e_{16}}$, $\br{e_2,e_7,e_{17}}$, $\br{e_2,e_8,e_{16}}$, $\br{e_3,e_7,e_{16}}$, $\br{e_7,e_7,e_{12}}$, and $\br{e_7,e_8,e_{11}}$ are all equal to 1.

By the Index Zero Axiom $\br{e_1,e_1,e_4}$, $\br{e_{1},e_{2},e_{3}}$ $\br{e_{2},e_{2},e_{2}}$, $\br{e_{8},e_{9},e_{9}}$ are equal to -2.

By the Composition Axiom $\br{e_7,e_9,ze_{10}}=\br{e_8,e_8,ze_{10}}=\pm 1$.

Consider the map $\varphi:\C[X,Y,Z]\rightarrow \h{S_{1,0}}{G}$ defined by $X\mapsto e_{16}$, $Y\mapsto e_1$ and $Z\mapsto e_7$ and extending to a $\C$-algebra homomorphism. One can check directly that this map is surjective. 

Our correlators tell us that $e_1e_7^4= 0$, $e_1^2=-2e_{16}$, $e_{16}e_1=e_{11}$,and that $e_7^5= -2e_{11}$. 

Hence $(2X+Y^2,Z^5+2XY,5YZ^4)\subset\ker\varphi$. Since $\C[X,Y,Z]/(2X+Y^2,Z^5+2XY,5YZ^4)=\q{S_{1,0}^T}$ has dimension 17, we deduce that the inclusion is equality, and so we have the isomorphism $\h{S_{1,0}}{G}\cong\q{S_{1,0}^T}$.
\newline\hrule

\subsubsection*{$\mathbf{S_{1,0}^T=x^2+yz^5+xy^2}$}
\[
\jac=(2x+y^2,z^5+2xy,5yz^4)
\phantom{XXX}
q_x=\tfrac {10}{20}, q_y=\tfrac 5{20}, q_z=\tfrac 3{20}, \hat c=\tfrac{24}{20}
\phantom{XXX}
G=\br J\cong \Z/20\Z
\]
\[
\Fix{J^k}=\begin{cases}
\C^3 & \text{if}\quad k=0\\
\C^2_{xy} & \text{if}\quad 4|k\\
\C_x  & \text{if}\quad 2|k, 4\nmid k\\
0    & \text{otherwise}
\end{cases}
\phantom{XXX}
\q{}|_{\Fix J^k}=
\begin{cases}
\br{1,y, y^2, z, z^2, \dots, z^8, yz,yz^2,yz^3,y^2z,y^2z^2,y^2z^3},\mu=17 \\
\br{1} \\
\br{1,y,y^2},\mu=3 \\
\br 1  \\
\end{cases}
\]

\begin{center}
\begin{tabular}{|c||c|c|c|c|c|c|c|c|c|c|c|c|c|c|c|c|c|c|c|c|c|c|c|c|}
\hline
$k$        & 1 & 3 & 4 & 5 & 7 & 8 & 9& 11 & 12 & 13 & 15 & 16 & 17 & 19\\
\hline
$|G|\cdot\deg_W$    & 0 & 32 & 28 & 24 & 16 & 12 & 8& 40 & 36 & 32 & 24 & 20 & 16 & 48\\
\hline
invariants& $\1$& $e_3$& $ye_4$& $e_5$& $e_7$& $ye_8$& $e_9$& $e_{11}$& $ye_{12}$& $e_{13}$& $e_{15}$& $ye_{16}$& $e_{17}$& $e_{19}$\\
\hline
\end{tabular}
\end{center}

Potential non-zero correlators:

By the pairing axiom $\br{\1,\1,e_{19}}$, $\br{\1,e_3,e_{17}}$,
$\br{\1,e_5,e_{15}}$, $\br{\1,e_7,e_{13}}$, and $\br{\1,e_9,e_{11}}$ are equal to 1, and $\br{ye_4,ye_{16},\1}=\br{ye_8,ye_{12},\1}=-1/2$.

By the concavity axiom $\br{e_3,e_9,e_9}$, $\br{e_5,e_7,e_9}$, $\br{e_7,e_{17},e_{17}}$, and $\br{e_9,e_{15},e_{17}}$ are equal to 1.

By the Index Zero axiom $\br{e_7,e_7,e_7}=-5$,

The composition axiom does not allow us to compute the last four three-point correlators, however if we let $a=\br{ye_4,e_8,e_9}$, $b=\br{e_5,ye_8,ye_8}$, $c=\br{ye_8,ye_{16},e_{17}}$, and $d=\br{ye_{16},ye_{16},e_9}$, we get the following relations:
\begin{eqnarray*}
b=-2ac\\
c=-2ad
\end{eqnarray*}
If $b\neq 0$, then we can construct the desired isomorphism. 

Consider the map $\varphi:\C[X,Y,Z]\rightarrow \h{S_{1,0}^T}{G}$ defined by $X\mapsto \beta y^3e_8$, where $\beta^2=\frac{-2}b$, $Y\mapsto e_9$ and $Z\mapsto e_7$ and extending to a $\C$-algebra homomorphism. One can check directly that this map is surjective. 

Our correlators tell us that $(y^3e_8)e_7= 0$, $e_7^2=-5e_{13}$, $e_9^4=e_{13}$, $(\beta y^3e_8)^2=(\beta^2b)e_{15}=-2e_{15}$ and that $e_9e_7=e_{15}$. 

Hence $(2XZ,Z^2+5Y^4,X^2+2YZ)\subset\ker\varphi$. Since $\C[X,Y,Z]/(2XZ,Z^2+5Y^4,X^2+2YZ)=\q{S_{1,0}}$ has dimension 14, we deduce that the inclusion is equality, and so we have the isomorphism $\h{S_{1,0}^T}{G}\cong\q{S_{1,0}}$ assuming $b\neq 0$.

\subsection{Singularities of Corank 2}\label{co2}
Here we consider only three singularities of the four singularities listed by Arnol'd, namely $J_{3,0}$, $Z_{1,0}$ and $W_{1,0}$. The singularity $W^{\#}_{1,2q}$ is not quasi-homogeneous for any choice of constants.

\subsubsection*{$\mathbf{J_{3,0}=x^3+bx^2y^3+y^9, 4b^3+27\neq 0}$}
\[
\jac=(3x^2+2bxy^3, 3bx^2y^2+9y^8)
\phantom{XXX}
q_x=\tfrac 3{9},\ \ q_y=\tfrac 19,\ \ \hat c=\tfrac{10}9
\phantom{XXX}
G=\br J\cong \Z/9\Z
\]

Notice that if $b=0$, then the maximal group of diagonal symmetries is not cyclic. But in that case $\h{Z_{3,0}}{G}|_{b=0}$ is isomorphic to a tensor product of simple singularities. So in this paper we only consider the case when the admissible group is generated by $J$.
\[
\Fix{J^k}=\begin{cases}
\C^2 & \text{if}\quad k=0\\
\C_x & \text{if}\quad 3|k\\
0    & \text{otherwise}
\end{cases}
\phantom{XXX}
\q{}|_{\Fix J^k}=
\begin{cases}
\br{1,x, y,y^2,\dots, y^7, xy,xy^2,\dots, xy^7},\mu=16 \\
\br{1,x},\mu=2 \\
\br1 
\end{cases}
\]

\begin{center}
\begin{tabular}{|c||c|c|c|c|c|c|c|}
\hline
$k$        & 0   & 1 & 2 & 4& 5& 7&8\\
\hline
$|G|\cdot\deg_W$         &10 &0 &8 &6 &14 &12 &20 \\
\hline
invariants& $y^5e_0,xy^2e_0$& $\1$& $e_2$& $e_4$& $e_5$& $e_7$& $e_8$\\
\hline
\end{tabular}
\end{center}

Potential non-zero correlators:

By of the pairing axiom $\br{\1\1,e_8}=\br{\1,e_2,e_7}=\br{\1,e_4,e_5}=1$, and
\begin{eqnarray*}
\br{y^5e_0,y^5e_0,\1}=\dfrac{2b^2}{9(4b^3+27)}\\
\br{y^5e_0,xy^2e_0,\1}=\dfrac{1}{4b^3+27}\\
\br{xy^2e_0,xy^2e_0,\1}=\dfrac{-2b}{3(4b^3+27)}
\end{eqnarray*}

By the concavity axiom $\br{e_2,e_4,e_4}=1$

Notice that in this case, we have the interesting situation that $b=0$ gives a different $FJRW$ ring for this singularity, than any other value for $b$.

Another interesting observation is that in this case, there is no Milnor ring that is isomorphic to the FJRW ring. To see this suppose there is a quasi-homogeneous polynomial $f(x_1,x_2,x_3,x_4)$ with 
\[
\C[x_1,x_2,x_3,x_4]/\jac_f\cong\h{J_{3,0}}{\br J} 
\]
and let $q_1,q_2,q_3,q_4$ be the corresponding charges. The isomorphism must send generators to generators, so we may assume $x_1\mapsto\mu$, $x_2\mapsto \nu$, $x_3\mapsto e_2$, and $x_4\mapsto e_4$, where $\mu$ and $\nu$ are in the untwisted sector.  Since we have an isomorphism of graded rings,
\[
q_1=10\alpha,\quad q_2=10\alpha,\quad q_3=8\alpha,\quad q_4=6\alpha,\quad \hat c_f=20\alpha
\]
(We admit the possibility that the grading may differ by a uniform scaling, although we see below that this possibility does not occur).

From the definition of $\hat c_f$, we have 
\[
\hat c_f=\sum(1-2q_i)=4-68\alpha,
\]
which we can solve to find $\alpha=\tfrac{1}{22}$, so the weights are 
\[
q_1=\tfrac{5}{11},\quad q_2=\tfrac{5}{11},\quad q_3=\tfrac{4}{11},\quad q_4=\tfrac{3}{11}
\]
Now this Milnor ring must have dimension $\mu=8$, but the dimension is given in terms of the charges by
\[
\mu=\prod(\tfrac{1}{q_i}-1)=\tfrac{168}{25},
\]
a contradiction. 
\newline\hrule

\subsection*{$\mathbf{Z_{1,0}=x^3y+y^7}$}
\[
\jac=(3x^2y,x^3+7y^6)
\phantom{XXX}
q_x = \tfrac 27,\ \ q_y = \tfrac 17,\ \ \hat c=\tfrac 87
\phantom{XXX}
G\cong \Z/21\Z
\phantom{XXX}
\br J\cong \Z/7\Z
\]

\[
\Fix \g^k=\begin{cases}
\C^2 & \text{if}\quad k=0\\
\C_y & \text{if}\quad 7|k, k\neq 0\\
0    & \text{otherwise}
\end{cases}
\phantom{XXX}
\q{}|_{\Fix \g^k}=
\begin{cases}
\br{1,x,x^2,y,\dots,y^6,xy,xy^2,\dots,xy^6}, \mu=21 \\
\br{1,y,y^2,y^3,y^4,y^5}, \mu=6  \\
\br 1 
\end{cases}
\]

\begin{center}
\begin{tabular}{|c||c|c|c|c|c|c|c|c|c|c|c|c|c|c|c|c|c|c|c|}
\hline
$k$        & 0  & 1    & 2    & 3    & 4   & 5   & 6& 8    & 9    & 10& 11  & 12 & 13  & 15  & 16   & 17 & 18 & 19 & 20\\
\hline
$|G|\cdot\deg_W$    &24 & 20& 16& 12& 8& 4& 0& 34& 30 &26& 22& 18& 14& 48& 44& 40& 36& 32& 28\\
\hline
invariants &$x^2e_0$ &$e_1$ &$e_2$ &$e_3$ & $e_4$ &$e_5$ &$\1$ &$e_8$ &$e_9$ &$e_{10}$& $e_{11}$ &$e_{12}$ &$e_{13}$ &$e_{15}$ &$e_{16}$ &$e_{17}$ &$e_{18}$ &$e_{19}$ &$e_{20}$ \\
\hline
\end{tabular}
\end{center}

Potential non-zero correlators:

By the pairing axiom $\br{\1,e_1,e_{20}}$, $\br{\1,e_2,e_{19}}$, $\br{\1,e_3,e_{18}}$, $\br{\1,e_4,e_{17}}$, $\br{\1,e_5,e_{16}}$, $\br{\1,\1,e_{15}}$, $\br{\1,e_8,e_{13}}$, $\br{\1,e_9,e_{12}}$, and $\br{\1,e_{10},e_{11}}$ are all equal to 1 and  $\br{x^2e_0, x^2e_0, \1}=-1/3$.

By the concavity axiom we have the following 3-point correlators are equal to 1: $\br{e_{1},e_{13},e_{13}}$, $\br{e_{2},e_{5},e_{20}}$, $\br{e_{2},e_{12},e_{13}}$, $\br{e_{3},e_{5},e_{19}}$, $\br{e_{3},e_{4},e_{20}}$, $\br{e_{3},e_{11},e_{13}}$, $\br{e_{3},e_{12},e_{12}}$, $\br{e_{4},e_{4},e_{19}}$, $\br{e_{4},e_{5},e_{18}}$, $\br{e_{4},e_{10},e_{13}}$, $\br{e_{4},e_{11},e_{12}}$, $\br{e_{5},e_{5},e_{17}}$, $\br{e_{5},e_{9},e_{13}}$, $\br{e_{5},e_{10},e_{12}}$, and $\br{e_{5},e_{11},e_{11}}$.

By the index zero axiom $\br{e_{1},e_{1},e_{4}}$, $\br{e_{1},e_{2},e_{3}}$, and $\br{e_{2},e_{2},e_{2}}$ are all equal to -3.

By the composition axiom using an argument similar that that for equation (\ref{e11}) we have
$\br{x^2e_{0},e_{1},e_{5}}=\br{x^2e_{0},e_{2},e_{4}}=\br{x^2e_{0},e_{3},e_{3}}= \pm 1$

Again, by examining degrees, we see that $e_5$ and $e_{13}$ are generators for $\h{J_{1,0}}{G}$.

Consider the map $\varphi:\C[X,Y]\rightarrow \h{J_{1,0}}{G}$ defined by $X\mapsto e_{13}$ and $Y\mapsto e_5$ and extending as a $\C$-algebra homomorphism. One can check directly that this map is surjective. 

Straightforward computations show that $e_5^7=-3e_{20}$, $e_{13}^2=e_{20}$, and $e_{13}e_5^6=0$. 

Hence $(Y^7+3X^2,XY^6)\subset \ker\varphi$. Since $\C[X,Y]/(Y^7+3X^2,XY^6)=\q{Z_{1,0}}$ has dimension 19, we deduce that the inclusion is equality. The reader will note that $Z_{1,0}$ is the transposed singularity for $E_{19}$ which we computed as an example in the introduction, so we have the isomorphism $\q{E_{19}}\cong\h{J_{1,0}}{G}$.
\newline\hrule

\subsubsection*{$\mathbf{W_{1,0}=x^4+ax^2y^3+y^6, a^2\neq 4}$}
\[
\jac=(4x^3+2axy^3, 3ax^2y^2+6y^5)
\phantom{XXX}
q_x=\tfrac 3{12},\ \ q_y=\tfrac 2{12},\ \ \hat c=\tfrac {14}{12}
\phantom{XXX}
G=\br J\cong \Z/12\Z
\]

Just as with $J_{3,0}$, if $a=0$, then the maximal group of diagonal symmetries is not cyclic. But in that case $\h{Z_{3,0}}{G}|_{b=0}$ is isomorphic to a tensor product of simple singularities. So in this paper we only consider the case when the admissible group is generated by $J$.
\[
\Fix{J^k}=\begin{cases}
\C^2 & \text{if} \quad k=0\\
C_x        & \text{if} \quad 4|k   \\
C_y        & \text{if} \quad 3|k  \\
0    & \text{otherwise}
\end{cases}
\phantom{XXX}
\q{}|_{\Fix J^k}=
\begin{cases}
\br{1, x, x^2, y,\dots,y^4,xy,\dots,xy^4,x^2y,\dots,x^2y^4},\mu=15 \\
\br{1,x,x^2},\mu=3 \\
\br{1,y,y^2,y^3,y^4},\mu=5 \\
\br 1  
\end{cases}
\]

\begin{center}
\begin{tabular}{|c||c|c|c|c|c|c|c|}
\hline
$k$        & 0   & 1 & 2 & 5 & 7 & 10 & 11 \\
\hline
$|G|\cdot\deg_W$    & 14  & 0 & 10 & 16 & 12 & 18 & 28\\
\hline
invariants& $xy^2e_0$& $\1$& $e_2$& $e_5$& $e_7$& $e_{10}$& $e_{11}$\\
\hline
\end{tabular}
\end{center}

Potential non-zero correlators:

By the Pairing axiom $\br{\1,\1,e_{11}}$, $\br{\1,e_{2},e_{10}}$, $\br{\1,e_{3},e_{9}}$, $\br{\1,e_{4},e_{8}}$, $\br{\1,e_{5},e_{7}}$, and $\br{\1,e_{6},e_{6}}$ are equal to 1 and $\br{xy^2e_{0},xy^2e_{0},\1}=\frac 1{24-6a^2}$

By the Concavity axiom $\br{e_{2},e_{2},e_{9}}$, $\br{e_{7},e_{9},e_{9}}$, and $\br{e_{8},e_{8},e_{9}}$ are equal to 1.

Just as with $J_{3,0}$, there is no Milnor ring isomorphic to this FJRW-ring. To see this, suppose there exists a quasi-homogeneous polynomial $f(x_1,x_2,\dots, x_5)$ with 
\[
\C[x_1,x_2,\dots,x_5]/\jac_f\cong\h{W_{1,0}}{\br J} 
\]
and let $q_1,q_2,\dots,q_5$ be the charges for $x_1,x_2,\dots,x_5$ resp. The isomorphism must send generators to generators, so we may assume $x_1\mapsto (xy^2e_0)$, $x_2\mapsto e_2$, $x_3\mapsto e_5$, $x_4\mapsto e_7$, and $x_5\mapsto e_{10}$. Since we have an isomorphism of graded rings (again allowing for a uniform rescaling of degrees), 
\[
q_1=14\alpha,\quad q_2=10\alpha,\quad q_3=16\alpha,\quad q_4=12\alpha,\quad q_5=18\alpha,\quad \hat c_f=28\alpha
\]
From the definition of $\hat c_f$, we have
\[
\hat c_f=\sum(1-2q_i)=5-140\alpha,
\]
Which we can solve to find $\alpha=\tfrac{5}{168}$, so the weights are 
\[
q_1=\tfrac{1}{12},\quad q_2=\tfrac{5}{84},\quad q_3=\tfrac{10}{21},\quad q_4=\tfrac{5}{14},\quad q_5=\tfrac{15}{28}.
\]
The dimension of the Milnor ring is given in terms of the charges by
\[
\mu=\prod(\tfrac{1}{q_i}-1),
\]
which is not even an integer.  So the FJRW A-model in this case is not isomorphic to the Milnor ring of a quasi-homogeneous polynomial.

\normalsize 
\bibliographystyle{acm}

\end{document}